\documentclass[11pt]{amsart}
\usepackage{amsmath,amssymb,amsfonts,url,mathptmx,tikz}
\usetikzlibrary{arrows}
\tikzset{ >=stealth'}
\numberwithin{equation}{section}
\newtheorem{Def}{Definition}[section]
\newtheorem{thm}{Theorem}[section]
\newtheorem{lem}{Lemma}[section]
\newtheorem{rem}{Remark}[section]
\newtheorem{prop}{Proposition}[section]

\begin{document}
\title[Toda system]{Classification of blowup limits for $SU(3)$ Singular Toda systems} 

\author{Chang-shou Lin}
\address{Department of Mathematics\\
        Taida Institute of Mathematical Sciences\\
        National Taiwan University\\
         Taipei 106, Taiwan } \email{cslin@math.ntu.edu.tw}

\author{Juncheng Wei}
\address{Department of Mathematics \\ University of British Columbia \\
 Vancouver, B.C., Canada V6T 1Z2 \\
 and Department of Mathematics\\
        Chinese University of Hong Kong\\
        Shatin, Hong Kong } \email{jcwei@math.ubc.ca}

\author{Lei Zhang}
\address{Department of Mathematics\\
        University of Florida\\
        358 Little Hall P.O.Box 118105\\
        Gainesville FL 32611-8105}
\email{leizhang@ufl.edu}

\date{\today}

\begin{abstract} For singular $SU(3)$ Toda systems, we prove that
 the limit of energy concentration is a finite set. In addition, for fully bubbling solutions we use Pohozaev identity to prove a
 uniform estimate. Our results extend previous results of Jost-Lin-Wang \cite{jostlinwang} on regular $SU(3)$ Toda systems.
\end{abstract}


\maketitle

\section{Introduction}

Systems of elliptic equations in two dimensional space with exponential nonlinearity are very commonly observed in Physics, Geometry, Chemistry and Biology. In this article we consider the following general system of equations defined in $\mathbb R^2$:
\begin{equation}\label{12aug14e1}
\Delta u_i+\sum_{j\in I}a_{ij}h_j e^{u_j}=4\pi \gamma_i\delta_0,\quad \mbox{ in }B_1\subset \mathbb R^2\quad \mbox{ for }i\in I,
\end{equation}
where $I=\{1,...,n\}$, $B_1$ is the unit ball in $\mathbb R^2$, $h_1,...,h_n$ are smooth functions, $A=(a_{ij})_{n\times n}$ is a constant matrix, $\gamma_i>-1$, $\delta_0$ is the Dirac mass at $0$.
If $n=1$ and $a_{11}=1$, the system (\ref{12aug14e1}) is reduced to a single Liouville equation, which has vast background in conformal geometry and Physics. The general system (\ref{12aug14e1}) is used for many models in different disciplines of science. If the coefficient matrix $A$ is non-negative, symmetric and irreducible, (\ref{12aug14e1}) is called a Liouville system and is related to models in the theory of chemotaxis \cite{child,keller}, in the Physics of charged particle beams \cite{bennet,debye,kiessling}, and in the theory of semi-conductors \cite{mock}, see
\cite{chanillo,chipot,linzhang1} and the reference therein for more applications of Liouville systems. If $A$ is the following Cartan matrix $A_n$:
$$A_n=\left(\begin{array}{ccccc}
2 & -1 & 0 & ... & 0\\
-1 & 2 & -1 & ... & 0\\
0 & -1 & 2 &  & 0\\
\vdots & \vdots & & \vdots\\
0 & \ldots & -1 & 2 & -1\\
0 & \ldots &  & -1 & 2
\end{array}
\right),
$$
the system (\ref{12aug14e1}) is called $SU(n+1)$ Toda system (which has $n$ equations) and is related to the non-abelian gauge in Chern-Simons theory, see \cite{dunne1,dunne2,ganoulis,leznov1,leznov2,MN,MR,mansfield,nolasco1,nolasco2,yang1,yang2} and the references therein. There are also many
works on the
relationship between $SU(n+1)$ Toda systems and holomorphic curves in $\mathbb{CP}^n$, flat $SU(n+1)$ connection, complete integrability and harmonic sequences, see \cite{bolton1,bolton2,calabi,chern,doliwa,guest,leznov2,lin-wei-ye} for references.

After decades of extensive study, many important questions related to the scalar  Liouville equation are answered and the behavior of blowup solutions is well understood (see \cite{bart2,bart3,bart4,chenlin1,chenlin2} for related discussions).  However, the understanding of blowup solutions to the more general systems (\ref{12aug14e1}) is far from complete.
In recent years, much progress has been made on more general systems and we only mention a few works related to the topic of the current article. First, Lin and Zhang \cite{linzhang1,linzhang2} completed a degree counting project for Liouville systems defined on Riemann surfaces. Second, for regular $SU(3)$ Toda systems (which have two equations), Jost-Lin-Wang \cite{jostlinwang} proved some uniform estimates for fully bubbling solutions (see section 4 for definition) using holonomy theory. Later Lin-Wei-Zhao \cite{lin-wei-zhao} improved the estimate of Jost-Lin-Wang to the sharp form using the non-degeneracy of the global $SU(3)$ solutions, which is established in Lin-Wei-Ye \cite{lin-wei-ye} among other things.

In this article we mainly focus on the asymptotic behavior of blowup solutions of (\ref{12aug14e1}) and the weak limit of energy concentration for $SU(n+1)$ Toda system.  More specifically,
let $u^k=(u_1^k,...,u_n^k)$ be a sequence solutions
\begin{equation}\label{12jun22e8}
\Delta u_i^k +\sum_{j=1}^n a_{ij} h_j^k e^{u_j^k}=4\pi \gamma_i^k\delta_0, \quad \mbox{ in }\quad B_1,\quad i=1,...,n
\end{equation}
with $0$ being its only possible blowup point in $B_1$:
\begin{equation}\label{12jun27e1}
\max_{K\subset\subset B_1\setminus \{0\}} u_i^k\le C(K).
\end{equation}
Since the right hand side of (\ref{12jun22e8}) is a Dirac mass, we define the regular part of $u_i^k$ to be
\begin{equation}\label{reg-u}
\tilde u_i^k(x)=u_i^k(x)-2\gamma_i^k \log |x|, \quad i=1,..,n, \quad x\in B_1.
\end{equation}
 $u^k=(u_1^k,...,u_n^k)$ is called \emph{a sequence of blowup solutions} if $\max_i\max_{x\in B_1}\tilde u_i^k\to \infty$.

We assume that $\gamma_i^k\to \gamma_i>-1$,
 $h_1^k,...,h_n^k$ are positive smooth functions with a uniform bound on their $C^3$ norm:
\begin{equation}\label{ah}
\frac 1C\le h_i^k\le C, \quad \|h_i^k\|_{C^3(B_1)}\le C, \quad \mbox{in } B_1, \quad \gamma_i^k\to \gamma_i>-1,\,\, \forall i\in I;
\end{equation}
and we suppose that there is a uniform bound on the oscillation of $u_i^k$  on $\partial B_1$ and its energy ( $\int_{B_1}h_i^ke^{u_i^k}$ is called the energy of $u_i^k$):
\begin{equation}\label{osc-enr}
|u_i^k(x)-u_i^k(y)|\le C, \quad \forall x,y\in \partial B_1,\qquad
\int_{B_1}h_i^k e^{u_i^k}\le C, \quad i\in I,
\end{equation}
where $C$ is independent of $k$.

Note that the oscillation finiteness assumption in (\ref{osc-enr}) is natural and generally satisfied in most applications. The energy bound in (\ref{osc-enr}) is also natural for systerm/equation defined in two dimensional space.

If $A=A_2$, system (\ref{12jun22e8}) describes $SU(3)$ with sources. Our first main theorem is concerned with the energy limits of solutions to
singular $SU(3)$ Toda systems.

Given any $\delta>0$, $u^k$ has no blowup point in $B_1\setminus B_{\delta}$ ( in this article we use $B(x,r)$ to denote a ball centered at $x$ with radius $r$ and use $B_r$ to denote $B(0,r)$). Thus we are interested in the following limit:
\begin{equation}\label{12sep2e3}
\sigma_i=\lim_{\delta\to 0}\lim_{k\to \infty} \frac{1}{2\pi} \int_{B_{\delta}} h_i^k e^{u_i^k},\quad i=1,2.
\end{equation}
Since for each $\delta>0$, $\int_{B_{\delta}}h_i^k e^{u_i^k}$ is uniformly bounded, the $\lim_{k\to \infty}$ in (\ref{12sep2e3}) is understood as the limit of a subsequence of $u^k$. For convenience we don't distinguish $u^k$ and its subsequences in this article.

Let
$$\mu_i=1+\gamma_i,\quad i=1,2 $$
and
$$\Gamma=\{(\sigma_1,\sigma_2):\,\, \sigma_1,\sigma_2\ge 0,\,\, \sigma_1^2-\sigma_1\sigma_2+\sigma_2^2=2\mu_1\sigma_1+2\mu_2\sigma_2 \} $$
be a quadratic curve in the first quadrant. It is easy to see that $\Gamma$ is contained in the box
$$[0,\frac 43 \mu_1+\frac 23\mu_2+\frac 43\sqrt{\mu_1^2+\mu_1\mu_2+\mu_2^2}]\times [0,\frac 23\mu_1+\frac 43\mu_2+\frac 43\sqrt{\mu_1^2+\mu_1\mu_2+\mu_2^2}]. $$
In Definition \ref{def1} below we shall define a finite set on $\Gamma$. In order to describe the mutual positions of points we say $(c,d)$ is \emph{in the upper right part of} $(a,b)$ if $c\ge a$ and $d\ge b$.

\begin{Def}\label{def1}
It is easy to verify that the following six points are on $\Gamma$:
\begin{eqnarray*}
(0,0),(2\mu_1,0),\,\, (0, 2\mu_2),\,\, (2\mu_1,2(\mu_1+\mu_2)),\\
(2(\mu_1+\mu_2),2\mu_2),\, \quad (2(\mu_1+\mu_2), 2(\mu_1+\mu_2)).
\end{eqnarray*}
First we let the six points above belong to $\Sigma$,
then we determine other points in $\Sigma$ as follows: For $(a,b)\in \Sigma$ intersect $\Gamma$ with $\sigma_1=a+2N$ and
$\sigma_2=b+2N$ ($N=0,1,2,...$) and add the point(s) of intersection to $\Sigma$ that belong to the upper right part of $(a,b)$. For each new member $(c,d)\in \Sigma$ added by this process, we apply the same procedure based on $(c,d)$ to obtain possible new members.
\end{Def}

\begin{thm}\label{cla-ener} Let $A=A_2$, $h_i^k$ and $\gamma_i^k$ satisfy (\ref{ah}). Then for $u^k$ satisfying (\ref{12jun22e8}), (\ref{12jun27e1}) and (\ref{osc-enr}), we have $(\sigma_1,\sigma_2)\in \Sigma$, where $\sigma_i$ is defined by (\ref{12sep2e3}) and $\Sigma$ is defined as in Definition \ref{def1}.
\end{thm}

\begin{rem} If $\gamma_1=\gamma_2=0$, the system is the nonsingular $SU(3)$ Toda system.  One sees easily that
$$\Sigma=\{(0,0),(2,0),(0,2),(2,4),(4,2),(4,4)\}. $$
Indeed, when the procedure described in Definition \ref{def1} is applied to any of the six points in $\Sigma$, no extra point of intersection can be found. For example if we start from $(0,0)$ and intersect $\Gamma$ by lines $\sigma_1=2N$ ($N$ being nonnegative integers). Then we see immediately that the intersection of $\Gamma$ with $\sigma_1=2$ gives $(2,0)$ and $(2,4)$, which are already in $\Sigma$. The intersection with $\sigma_1=4$ gives $(4,2)$ and $(4,4)$, which also belong to the six types in $\Sigma$. There is no intersection between $\Gamma$ and $\sigma_1=6$.
Theorem \ref{cla-ener} in this special case was proved by Jost-Lin-Wang in \cite{jostlinwang}.
A recent work of Pistoia-Musso-Wei \cite{pmw} proved that all the six cases for nonsingular $SU(3)$ Toda systems can occur.
\end{rem}

		
	
	
					


\begin{rem} It is easy to observe that the maximum value of $\sigma_1$ on $\Gamma$ is
$$\frac 43\mu_1+\frac 23\mu_2+\frac 43\sqrt{\mu_1^2+\mu_1\mu_2+\mu_2^2}. $$ The maximum value of $\sigma_2$ is
$$\frac 23\mu_1+\frac 43\mu_2+\frac 43\sqrt{\mu_1^2+\mu_1\mu_2+\mu_2^2}. $$
 Thus $\Sigma$ is a finite set. As two special cases, we see that
\begin{enumerate}
\item If
\begin{eqnarray*}
\frac 43\mu_1+\frac 23\mu_2+\frac 43\sqrt{\mu_1^2+\mu_1\mu_2+\mu_2^2}<2, \quad \mbox{ and }\\
\frac 23\mu_1+\frac 43\mu_2+\frac 43\sqrt{\mu_1^2+\mu_1\mu_2+\mu_2^2}<2.
\end{eqnarray*}
then there are only six points in $\Sigma$:
\begin{eqnarray*}
\Sigma&=&\Big\{(0,0),(2\mu_1,0),(0,2\mu_2),(2(\mu_1+\mu_2),2\mu_2),\\
&&(2\mu_1,2(\mu_1+\mu_2)),(2(\mu_1+\mu_2),2(\mu_1+\mu_2)) \Big\}.
\end{eqnarray*}
\item
For $\gamma_1=\gamma_2=1$, in addition to
$(0,0),(4,0),(0,4),(4,8),(8,4),(8,8)$, $\Sigma$ has other 14 points.
\end{enumerate}
\end{rem}

An earlier version of the current article was posted on arxiv.org in March 2013. After that some works have been done based on Theorem \ref{cla-ener} (see \cite{malchiodi-b} for example). Theorem \ref{cla-ener} reflects some essential differences between Toda systems and Lioiuville systems. Lin-Wei-Ye \cite{lin-wei-ye} proved that all the global solutions of $SU(n+1)$ Toda systems can be described by $n^2+2n$ parameters and the energy of global solutions is a discrete set. On the other hand, the global solutions of Liouville systems all belong to a family of three parameters but their energy forms a $n-1$ dimensional hyper-surface (see \cite{chipot,linzhang1}). These differences lead to very different approaches in their respective research. For example, Lin-Wei-Zhao \cite{lin-wei-zhao} obtained sharp estimates for fully-bubbling solutions (see section 4 for definition) of $SU(3)$ Toda system using the discreteness of energy as a key ingredient in their proof.

Here we briefly describe the strategy to prove Theorem \ref{cla-ener}. First we introduce a selection process suitable for $SU(n+1)$ Toda systems. The selection process has been widely used for prescribing curvature type equations (see \cite{li}, \cite{lin}, etc) and we modify it to locate the bubbling area, which is a union of finite disks. In each of the disks, the blowup solutions have roughly the energy of a global $SU(m+1)$ Toda system on $\mathbb R^2$ (with $m\le n$), which is the limit of the blowup solutions after scaling. If $m=n$, which means no component is lost after scaling and taking the limit, we say the sequence of solutions in the disk is fully bubbling. Otherwise we call it partially bubbling. Next we introduce the concept ``group" to place bubbling disks according to their mutual locations. There are only finite bubbling disks and their mutual distances may tend to $0$ with very different speed. The name ``group" is used to describe a few disks that are ro!
 ughly closest to one another and much further from other disks. Lemma \ref{oscillation} is a Harnack type result that plays an important role in determining the energy concentration around a group. Suppose there is a circle that surrounds a group and both components of the blowup solutions have fast decay (see section three for definition) on the circle. Then a Pohozaev identity can be computed on this circle to determine how much energy this group carries. Because of Lemma \ref{oscillation} such a circle can always be found, so the energy within the circle can be determined. Then we consider the combination of groups by scaling. The relationship among groups is similar to that of members in a same group. For example if the distance between two groups is scaled to be $1$, the bubbling disks of one group look like a Dirac mass from afar. We can similarly find circles surrounding groups that are also suitable for computing Pohozaev identities ( i.e. both components of the blo!
 wup solutions have fast decay on these circles). From these Pohozaev identities we determine how much energy is contained in each group and all the combinations of groups. One important fact is that one component of the blowup solutions always has fast decay, even though the second component may not be the case. It is possible for the first (fast decay) component to turn to a slow decay component as the distance to a group becomes bigger, but before that happens the second component, which used to be a slow decay component, will turn to fast decay component first.

As another application of the Pohozaev identity we establish some uniform estimates for fully bubbling solutions.
These estimates  were first obtained by Li \cite{licmp} for scalar Liouville equation without singularity (using the method of moving planes) and Bartolucci-et al. \cite{bclt} for scalar Liouville equation with singularity (using Pohozaev identity and potential analysis). For regular $SU(3)$ Toda systems Jost-Lin-Wang \cite{jostlinwang} established similar estimates using holonomy theory. Our results (Theorem \ref{thm3} and Theorem \ref{thm2}) apply to general $SU(n+1)$ Toda systems with singularity.

This article is set out as follows. In section two we introduce the selection process mentioned before and in section three we prove the Pohozaev identity, which is crucial for the proof of Theorem \ref{cla-ener}. Then in Section four we prove a uniform estimate for fully bubbling solutions (Theorem \ref{thm2} and Theorem \ref{thm3}). Then in section five and section six we finish the proof of Theorem \ref{cla-ener} according to the strategy mentioned before.

\medskip

{\bf Acknowledgements:} Part of the paper was finished when the third author was visiting Chinese University of Hong Kong in April, May and December in 2012, and Taida Institute of Mathematical Sciences (TIMS) in June 2012. He would like to thank both institutes for their warm hospitality and generous financial support. The research of J. Wei is partially supported by NSERC of Canada.

\section{A selection process for $SU(n+1)$ Toda systems}

Clearly in the proof of Theorem \ref{cla-ener} we can assume $0$ to be a blowup point:
\begin{equation}\label{did-bu}
\max_{x\in B_1,i\in I} \{u_i^k-2\gamma_i^k \log |x| \}\to \infty
\end{equation}
because otherwise the blowup type is $(0,0)$. So from now on throughout the paper (\ref{did-bu}) is assumed.

\subsection{ Case one: $\gamma_1^k=...=\gamma_n^k=0$}

\begin{prop}\label{selectprop}
Let $A=(a_{ij})_{n\times n}$ be the Cartan matrix $A_n$, $h_i^k$ satisfy (\ref{ah}) and $u^k=(u_1^k,..,u_n^k)$ be a sequence of solutions to (\ref{12jun22e8}) with $\gamma_1^k=..=\gamma_n^k=0$
such that (\ref{osc-enr}) and (\ref{12jun27e1}) hold. Then there exist finite sequences of points $\Sigma_k:=\{x_1^k,....,x_m^k\}$ (all $x_j^k\to 0, j=1,...,m$) and positive numbers $l_1^k,...,l_m^k\to 0$ such that the following four properties hold:
\begin{enumerate}\item
$\max_{i\in I}\{ u_i^k(x_j^k)\}=\max_{B(x_j^k,l_j^k),i\in I}\{ u_i^k\}$ for all $j=1,..,m$.
\item
$exp (\frac 12 \max_{i\in I}\{u_i^k(x_j^k)\}) l_j^k\to \infty, \quad j=1,...,m$.
\item There exists $C_1>0$ independent of $k$ such that
$$ u_i^k(x)+2\log \,\, dist(x,\Sigma_k)\le C_1, \quad \forall x\in B_1,\quad i\in I $$
where $dist$ stands for distance.
\item In each $B(x_j^k,l_j^k)$ let
\begin{equation}\label{vik-sel}
v_i^k(y)=u_i^k(\epsilon_k y+x_j^k)+2\log \epsilon_k, \quad \epsilon_k=e^{-\frac 12 M_k},\quad M_k=\max_i\max_{B(x_j^k,l_j^k)}u_i^k.
\end{equation}
 Then one of the following two alternatives holds\\
 (a):\quad The sequence is fully bubbling: along a subsequence $(v_1^k,...,v_n^k)$ converges in $C^2_{loc}(\mathbb R^2)$ to $(v_1,...,v_n)$ which satisfies
$$\Delta v_i+\sum_{j\in I} a_{ij}h_j e^{v_j}=0,\quad \mathbb R^2, \quad i\in I. $$
$$\lim_{k\to \infty} \int_{B(x_j^k,l_j^k)} \sum_{t\in I} a_{it}h_t^ke^{u_t^k}>4\pi,\quad i\in I.
$$
  (b):$I=J_1\cup J_2\cup... \cup J_m\cup N$ where $J_1,J_2,...,J_m$ and $N$ are disjoint sets. $N\neq \emptyset$ and each $J_t$ ($t=1,..,m$) consists of consecutive indices.  For each $i\in N$,  $v_j^k$ tends to $-\infty$ over any fixed compact subset of $\mathbb R^2$. The components of $v^k=(v_1^k,...,v_n^k)$ corresponding to each $J_l$ ($l=1,...,m$) converge in $C^2_{loc}(\mathbb R^2)$ to a $SU(|J_l|+1)$ Toda system, where $|J_l|$ is the number of indices in $J_l$. For each $i\in J_l$, we have
  $$\lim_{k\to \infty}\int_{B(x_j^k,l_j^k)}\sum_{t\in J_l}a_{it}h_t^ke^{v_t^k}>4 \pi. $$
\end{enumerate}
\end{prop}

\begin{rem} In this article we don't use different notations for sequences and subsequences.
\end{rem}

\begin{rem}\label{dis} For each $x_j^k\in \Sigma_k$ suppose $2t_j^k$ is the distance from $x_j^k$ to $\Sigma_k\setminus \{x_j^k\}$. Then $t_j^k/l_j^k\to \infty$
as $k\to \infty$ if $l_j^k$ is suitably chosen.
\end{rem}

\noindent{\bf Proof of Proposition \ref{selectprop}:}

Without loss of generality we assume
$$u_1^k(x_1^k)=\max_{i\in I,x\in B_1}u_i^k(x). $$
Clearly $x_1^k\to 0$ because $\max_i \max_{x\in B_1}u^k_i\to \infty$ and $u^k$ is uniformly bounded above away from the origin.
Let $( v_1^k,..,v_n^k)$ be defined by (\ref{vik-sel}) with $x_j^k$ replaced by $x_1^k$.
Immediately we observe that $|\Delta v_i^k|$ is bounded because each $v_i^k\le 0$. Consequently
 $|v_i^k(z)-v_i^k(0)|$
is uniformly bounded in any compact subset of $\mathbb R^2$. Thus, by $v_1^k(0)=0$, at least (along a subsequence) $v_1^k$ converges in $C^2_{loc}(\mathbb R^2)$ to a function $v_1$. For other components of $v^k=(v_1^k,.., v_n^k)$, either some of them tend to $-\infty$ over any compact subset of $\mathbb R^2$, or all of them converge to a system of $n$ equations.  Let $J\subset I$ be the set of indices corresponding to those convergent components. That is, for all $i\in J$, $v_i^k$ converges to $v_i$ in $C^2_{loc}(\mathbb R^2)$ and for all $j\in I\setminus J$, $v_i^k$ tends to $-\infty$ over any fixed compact subset of $\mathbb R^2$.
For each $i\in I\setminus J$, there is $J_1\subset J$ such that $i\in J_1$, the indices in $J_1$ are consecutive and the limit of $v_i^k$ is one component of a
$SU(|J_1|+1)$ Toda system:
\begin{equation}\label{12sep2e1}
\left\{\begin{array}{ll}
\Delta v_m+\sum_{j\in J}a_{ml}h_l e^{v_l}=0,\quad \mbox{ in } \mathbb R^2, \quad \forall m\in J_1 \\
\\
\int_{\mathbb R^2}h_me^{v_m}\le C,\quad m\in J_1
\end{array}
\right.
\end{equation}
where $h_m=\lim_{k\to \infty}h_m^k(x_1^k)$, $(a_{ij})=A_{|J_1|}$, and $C$ is the same constant as in (\ref{osc-enr}).
By the classification theorem of Lin-Wei-Ye \cite{lin-wei-ye} (if the limit is a system) or Chen-Li \cite{chenliduke} (if the limit is one equation) we have
\begin{equation}\label{12aug3e1}
\sum_{j\in J_1}\int_{\mathbb R^2} a_{ij} h_j e^{v_j}=8\pi,\quad \forall i\in J_1
\end{equation}
and
\begin{equation}\label{vi-loc}
v_i(x)=-4\log |x|+O(1),\quad |x|>2,\quad \forall i\in J_1.
\end{equation}
Thus for any index $i\in I$ we can find $R_k\to \infty$ such that
\begin{equation}\label{12sep2e2}
 v_i^k(y)+2\log |y|\le C, \quad |y|\le R_k,\quad \mbox{ for } i\in I.
\end{equation}
Equivalently for $u^k$ there exist $l_1^k\to 0$ such that
$$u_i^k(x)+2\log |x-x_1^k|\le C, \quad |x-x_1^k|\le l_1^k, \quad \mbox{for }\,\, i\in I $$
and
$$e^{\frac 12u_1^k(x_1^k)}l_1^k\to \infty,\quad i\in J,\quad  \mbox{ as }\,\, k\to \infty. $$

Next we let $q_k$ be the maximum point of $\max_{|x|<1,i\in I}u_i^k(x)+2\log |x-x_1^k|$.
If
$$\max_{|x|\le 1, i\in I}u_i^k(x)+2\log |x-x_1^k| \to \infty, $$  we let $j$ be the index such that
$$u_j^k(q_k)+2\log |q_k-x_1^k|=\max_{i\in I}u_i^k(x)+2\log |x-x_1^k|\to \infty. $$
The following localization is to adapt the original argument of R. Schoen \cite{schoen} for the scalar curvature equation (also see \cite{chen-lin-2}).
Set
$$d_k=\frac 12|q_k-x_1^k| $$
and
$$S_i^k(x)=u_i^k(x)+2\log \bigg (d_k-|x-q_k| \bigg )\,\, \mbox{in }\,\, B(q_k,d_k). $$
Then clearly for fixed $k$, $S_i^k\to -\infty$ as $x$ tends to $\partial B(q_k,d_k)$. On the other hand, at least for $j$ we have
$$S_j^k(q_k)=u_j^k(q_k)+2\log d_k\to \infty. $$
Let $p_k$ be where
$$\max_i \max_{x\in \bar B(q_k,d_k)}S_i^k  $$
is attained and $i_0$ be the index corresponding to where the the maximum is taken:
\begin{equation}\label{12aug17e1}
u_{i_0}^k(p_k)+2\log \bigg (d_k-|p_k-q_k| \bigg )\ge S_j^k(q_k)\to \infty.
\end{equation}
Let
$$l_k=\frac 12(d_k-| p_k -q_k |). $$
Then for $y\in B(p_k,l_k)$, by the choice of $p_k$ and $l_k$ we have
$$u_i^k(y)+2\log (d_k-|y-q_k|)\le u_{i_0}^k(p_k)+2\log (2l_k),\quad \forall i\in I. $$
On the other hand, by the definition of $l_k$ we have
$$d_k-|y-q_k|\ge d_k-|p_k-q_k|-|y-p_k|\ge l_k,\quad \mbox{ if }\,\, |y-p_k|<l_k, $$ and
\begin{equation}\label{star-1}
u_i^k(y)\le u_{i_0}^k(p_k)+2\log 2, \quad \forall y\in B(p_k,l_k).
\end{equation}
Next we set
\begin{equation}\label{star-2}
\mathcal{R}_k=e^{\frac 12 u_{i_0}^k(p_k)}l_k
\end{equation}
and scale $u_i^k$ by
$$\tilde v_i^k(y)=u_i^k(p_k+e^{-\frac 12u_{i_0}^k(p_k)}y)-u_{i_0}^k(p_k),\quad \mbox{ for } i\in I. $$
From (\ref{12aug17e1}) we clearly have $\mathcal{R}_k\to \infty$. By (\ref{star-1}) and standard elliptic estimates for the Laplace, $\tilde v_i^k$ is bounded in $C^2_{loc}(\mathbb R^2)$ and there exists
$\emptyset \neq J\subset I$ such that for all $i\in J$, $\tilde v_i^k$ converges to a limit system like (\ref{12sep2e1}). On the other hand $\tilde v_i^k$ converges uniformly to $-\infty$ over all compact subsets of $\mathbb R^2$ for all $i\in I\setminus J$. Clearly (\ref{12sep2e2}) holds for $\tilde v_i^k$. Going back to $u^k$ we have
$$u_i^k(x)+2\log |x-x_2^k|\le C, \mbox{ for } |x-x_2^k|\le l_2^k $$
where $x_2^k$ is where $max_i\max_{B(p_k,l_2^k)}u_i^k$ is attained and $l_2^k=l_k$. Here we note that $x_2^k$ is neither $q_k$ nor $p_k$ and the distance between $p_k$ and $x_2^k$ is small: $e^{\frac 12 u_{i_0}^k(p_k)}|x_2^k-p_k|=O(1)$.
If we re-scale $u^k$ around $x_2^k$, $v^k$ defined as in (\ref{vik-sel}) satisfies (a) and (b) in Proposition \ref{selectprop}. Clearly
$B(x_1^k,l_1^k)\cap
B(x_2^k,l_2^k)=\emptyset$.

To continue with the selection process, we let $\Sigma_{k,2}:=\{x_1^k,x_2^k\}$ and consider
$$\max_{i\in I,x\in B_1} u_i^k(x)+2\log dist (x, \Sigma_{k,2}). $$
If along a subsequence, the quantity above tends to infinity we apply the same procedure to get $x_3^k$ and $l_3^k$. Since after each selection we add a new disjoint disk, say $B(x_m^k,l_m^k)$, in which the profile of bubbling solutions is like that of a global system. From (\ref{12aug3e1}) we see that
 $$\int_{B(x_m^k,l_m^k)}\sum_i h_i^k e^{u_i^k}\ge C,  \quad \mbox{ for some } C>0 \mbox{ independent of } k. $$
 Therefore the process stops after finite steps by (\ref{osc-enr}).  Eventually we let
$$\Sigma_k=\{ x_1^k,...,x_L^k \} $$
and it holds
\begin{equation}\label{harnack1}
u_i^k(x)+2\log d(x, \Sigma_k)\le C, \quad i\in I.
\end{equation}
Proposition \ref{selectprop} is established. $\Box$

\medskip

\subsection{ Case two: Singular case $\exists \gamma_i\neq 0$. \quad }

First the selection process is almost the same. The difference is instead of taking the maximum of $u_i^k$ over $B_1$ we let $0\in \Sigma_k$. Clearly in $B_1\setminus \{0\}$
 $u^k$ satisfies the same equation as the nonsingular case.
 Then we consider the maximum of $u_i^k(x)+2\log dist(x, \Sigma_k)=u_i^k(x)+2\log |x|$ and the selection proceeds the same as before. Therefore
in the singular case $\Sigma_k=\{0, x_1^k,...,x_m^k\}$.

\begin{lem}\label{oscillation} Let $\Sigma_k$ be the blowup set (Thus if all $\gamma_i^k=0$, $\Sigma_k=\{x_1^k,...,x_m^k\}$, if the system is singular,
$\Sigma_k=\{0,x_1^k,...,x_m^k\}$). In either case for all $x_0\in B_1\setminus \Sigma_k$, there exists $C_0$ independent of $x_0$ and $k$ such that
$$|u_i^k(x_1)-u_i^k(x_2)|\le C_0,\quad \forall x_1,x_2\in B(x_0,d(x_0,\Sigma_k)/2),\mbox{ for all } i\in I. $$
\end{lem}

\noindent{\bf Proof of Lemma \ref{oscillation}:} We can assume $|x|<\frac 1{10}$ because it is easy to see from the Green's representation formula that the oscillation of $u_i^k$ on $B_1\setminus B_{1/10}$ is finite.
Recall the regular part of $u_i^k$ is defined in (\ref{reg-u}) and $\tilde u_i^k$ satisfies
$$\Delta \tilde u_i^k(x)+\sum_{j\in I}a_{ij}h_j^k(x)|x|^{2\gamma_j^k}e^{\tilde u_j^k(x)}=0, \quad B_1, \quad i\in I. $$
Let $\sigma_k$ be the distance between $x_0$ and $\Sigma_k$. Clearly, for $x_0\in B_1\setminus \Sigma_k$ and $x_1,x_2\in B(x_0,d(x_0,\Sigma_k)/2)$,
\begin{eqnarray*}
&&u_i^k(x_1)-u_i^k(x_2)\\
&=&\tilde u_i^k(x_1)-\tilde u_i^k(x_2)+O(1)\\
&=&\int_{B_1}(G(x_1,\eta)-G(x_2,\eta))\sum_{j\in I}a_{ij}h_j^k(\eta)|\eta|^{2\gamma_j^k}e^{\tilde u_j^k(\eta)}d\eta+O(1).
\end{eqnarray*}
Here $G$ is the Green's function on $B_1$. The last term on the above is $O(1)$ because it is the difference of two points of a harmonic function that has bounded oscillation on $\partial B_1$. Since both $x_1,x_2\in B_{1/10}$, it is easy to use the uniform bound on the energy (\ref{osc-enr}) to obtain
$$\int_{B_1}(\gamma(x_1,\eta)-\gamma(x_2,\eta))\sum_{j\in I}a_{ij}h_j^k(\eta)|\eta|^{2\gamma_j^k}e^{\tilde u_j^k(\eta)}d\eta=O(1)$$
where $\gamma(\cdot,\cdot)$ the regular part of $G$. Therefore we only need to show
$$\int_{B_1}\log \frac{|x_1-\eta|}{|x_2-\eta|}\sum_j a_{ij}h_j^k|\eta|^{2\gamma_j}e^{\tilde u_j}d\eta=O(1). $$
If $\eta\in B_1\setminus B(x_0, \frac 34\sigma_k)$, we have $\log (|x_1-\eta|/|x_2-\eta|)=O(1)$, then the integration over $B_1\setminus B(x_0,
\frac 34\sigma_k)$ is uniformly bounded. Therefore we only need to show
\begin{eqnarray*}
&&\int_{B(x_0,\frac 34\sigma_k)}\log \frac{|x_1-\eta|}{|x_2-\eta|}\sum_j a_{ij} h_j^k|\eta|^{2\gamma_j}e^{\tilde u_j^k}d\eta\\
&=&\int_{B(x_0,\frac 34\sigma_k)}\log \frac{|x_1-\eta|}{|x_2-\eta|}\sum_j a_{ij} h_j^ke^{u_j^k}d\eta=O(1).
\end{eqnarray*}
To this end, let
\begin{equation}\label{121222e1}
v_i^k(y)=u_i^k(x_0+\sigma_ky)+2\log \sigma_k, \quad i\in I,\quad y\in B_{3/4}.
\end{equation}
Then we just need to show
\begin{equation}\label{12may11e2}
\int_{B_{3/4}}\log \frac{|y_1-\eta|}{|y_2-\eta |} \sum_j a_{ij}h_j^k(x_0+\sigma_k\eta)e^{v_j^k(\eta)}d\eta=O(1).
\end{equation}
We assume, without loss of generality that
$e_1$ is the image of the closest blowup point in $\Sigma_k$. Thus by the selection process
$$ v_i^k(\eta)\le -2\log |\eta -e_1|+C. $$
Therefore
$$e^{ v_i^k(\eta)}\le C|\eta -e_1|^{-2}. $$
With this estimate we observe that
$|\eta -e_1|\ge C>0$ for $\eta\in B_{3/4}$. Thus for $j=1,2$ and any fixed $i\in I$,
$$\int_{B_{3/4}}\bigg | \log |y_j-\eta | \bigg |e^{v_i^k(\eta )}d\eta\le C\int_{B_{3/4}}\frac{\big |\log |y_j-\eta | \big |}{|\eta -e_1|^2}d\eta
\le C.$$
Lemma \ref{oscillation} is established. $\Box$

\begin{rem} For systems with nonnegative coefficient matrix $A$, the selection process can also be applied. See Chen-Li \cite{chen-li-duke2} or Lin-Zhang \cite{linzhang1} for more details.
\end{rem}

\section{Pohozaev identity and related estimates on the energy}

In this section we derive a Pohozaev identity for $u^k$ satisfying (\ref{12jun22e8}), (\ref{12jun27e1}), (\ref{osc-enr}),  $h_i^k$ and $\gamma_i^k$ satisfying (\ref{ah}), and $A=A_n$.

\begin{prop}\label{piforU} Let $A=A_n$, $\sigma_i$ be defined by (\ref{12sep2e3}).  Suppose $u^k=(u_1^k,...,u_n^k)$ satisfy
(\ref{12jun22e8}), (\ref{osc-enr}),(\ref{12jun27e1}) and (\ref{did-bu}), $h^k$ and $\gamma_i^k$ satisfy (\ref{ah}).
Then we have
$$\sum_{i,j\in I} a_{ij} \sigma_i\sigma_j=4\sum_{i=1}^n (1+\gamma_i)\sigma_i. $$
\end{prop}

\noindent{\bf Proof of Proposition \ref{piforU}:}

\begin{lem} \label{12aug8lem1}
Given any $\epsilon_k\to 0$ such that $\Sigma_k\subset B(0,\epsilon_k/2)$,
there exist $l_k \to 0$ satisfying $l_k\ge 2\epsilon_k$ and
\begin{equation}\label{12may23e1}
\bar u_i^k(l_k)+2\log l_k\to -\infty, \mbox{ for all } i\in I, \mbox{ where } \bar u_i^k(r):=\frac 1{2\pi r}\int_{\partial B_r}u_i^k.
\end{equation}
\end{lem}

\begin{rem} By Lemma \ref{12aug8lem1} and Lemma \ref{oscillation}
$$u_i^k(x)+2\log |x|\to -\infty,\quad  \forall i\in I \mbox{ and } \quad \forall x\in \partial B_{l_k}. $$
This is crucial for evaluating the $\mathcal{R}_1$ term ( the first term on the right) of (\ref{poho1}) below.
\end{rem}

\noindent{\bf Proof of Lemma \ref{12aug8lem1}:}
  Since
 $\Sigma_k\subset B(0, \epsilon_k/2)$, we have, by the third statement of Proposition \ref{selectprop},
 \begin{equation}\label{control1}
u_i^k(x)+2\log |x|\le C,\quad |x|\ge \epsilon_k.
\end{equation}
The key point of the argument below is that we can always use the finite energy assumption and Lemma \ref{oscillation} to make $u_1^k$ satisfy (\ref{12may23e1}). Then we can adjust the radius to make other components satisfy (\ref{12may23e1}) as well.

First we observe that for each fixed $i$ there exists $r_{k,i}\ge \epsilon_{k}$ such that
 \begin{equation}\label{12aug17e2}
\bar u_i^k(r_{k,i})+2\log r_{k,i}\to -\infty,
\end{equation}
because otherwise we would have
$$\bar u_i^k(r)+2\log r\ge -C  \quad \mbox{ for all } r\ge \epsilon_k$$
for some $C>0$. By Lemma \ref{oscillation} $u_i^k$ has bounded oscillation on each $\partial B_r$. Thus
$$u_i^k(x)+2\log |x|\ge -C \quad \mbox{ for all } x\in \partial B_r,\quad \epsilon_k<r< 1$$
for some $C$. Then
$$e^{u_{i}^k(x)}\ge C|x|^{-2},\quad \epsilon_{k}\le |x|\le 1. $$
Integrating $e^{u_{i}^k}$ on $B_1\setminus B_{\epsilon_{k}}$ we get a contradiction on the uniform energy bound of $\int_{B_1}h_i^k e^{u_i^k}$.
(\ref{12aug17e2}) is established.

 \medskip

First for $u_1^k$, we find $r_{k,1}\ge \epsilon_k$ so that
$$\bar u_1^k(r_{k,1})+2\log r_{k,1}\to -\infty. $$
Here we claim that we can assume $r_{k,1}\to 0$ as well. In fact, if $r_{k,1}$ does not tend to $0$, by Lemma \ref{oscillation}
$$\bar u_1^k(r)+2\log r\le -N_k+C, \quad r_{k,1}/2<r<r_{k,1} $$
where $N_k\to \infty$ and satisfies
$$\bar u_1^k(r_{k,1})+2\log r_{k,1}\le -N_k. $$
Using Lemma \ref{oscillation} again we have
$$\bar u_1^k(r)+2\log r\le -N_k+C, \quad r_{k,1}/4<r<r_{k,1}/2. $$
Obviously this process can be done $\bar N_k$ times where $\bar N_k$ is chosen to tend to infinity slowly enough so that
$\bar r_k:=r_{k,1}2^{-\bar N_k}$ satisfies
$$\bar u_1^k(\bar r_k)+2\log \bar r_k \le -N_k+C \bar N_k\to -\infty. $$
We can use $\bar r_k$ to replace $r_{k,1}$.  Exactly the same argument clearly shows the existence of $s_k\to 0$, $\tilde N_k\to \infty$ such that
$$\left\{\begin{array}{ll} s_k/r_{k,1}\to \infty,\\
\bar u_1^k(r)+2\log r\le -\tilde N_k,\quad r_{k,1}\le r\le s_k.
\end{array}
\right.
$$
Next we claim that between $r_{k,1}$ and $s_k$, there must be a $r_{k,2}$ such that
\begin{equation}\label{u2k-small}
\bar u_2^k(r_{k,2})+2\log r_{k,2}\le -N_{k,2}
\end{equation}
for some $N_{k,2}\to \infty$ as $k\to \infty$. The proof of (\ref{u2k-small}) is very similar to what has been used before: If this is not the case, $e^{u_2^k}\ge Cr^{-2}$ for some $C>0$ and $r\in (r_{k,1},s_k)$. The fact that $s_k/r_{k,1}\to \infty$ leads to a contradiction to the uniform bound of $u_2^k$'s energy.

Thus we have proved that for $r=r_{k,2}$ both $u_1^k,u_2^k$ decay faster than $-2\log r$:
$$\bar u_i^k(r)+2\log r\le -N_k,\quad i=1,2,\quad r=r_{k,2} $$
for some $N_k\to \infty$. Then it is easy to see that there exists $s_k\to 0$ and $s_k/r_{k,2}\to \infty$ such that
$$\bar u_i^k(r)+2\log r\le -N_k',\quad i=1,2,\quad r_{k,2}\le r\le s_k  $$
for some $N_k'\to \infty$ as well. The same argument above guarantees the existence of $l_k\in (r_{k,2},s_k)$ and some $N_k''\to \infty$ such that
$$\bar u_3^k(l_k)+2\log l_k\le -N_k''. $$
Clearly this argument can be applied finite times to exhaust all the components of the whole system.
Lemma \ref{12aug8lem1} is established. $\Box$

\medskip

Now we continue with the proof of Proposition \ref{piforU}.

\medskip

{\bf Case one: $\gamma_i^k\equiv 0$.}

\medskip

Using the definition of $\sigma_i$ in (\ref{12sep2e3}) we
choose $l_k\to 0$ such that $\Sigma_k\subset B(0, l_k/2)$  and
\begin{equation}\label{13mar10e1}
\frac{1}{2\pi }\int_{B_{l_k}}h_i^k e^{u_i^k}=\sigma_i+o(1), \quad \mbox{ for } i\in I.
\end{equation}
Here we claim that (\ref{12may23e1}) also holds, because otherwise we would have
$$\bar u_i(l_k)+2\log l_k\ge -C. $$
By Lemma \ref{oscillation}
$$\bar u_i(r)+2\log r\ge -C_1,\quad l_k\le r\le 2l_k, $$
which means there is a lower bound on the energy in the annulus $B_{2l_k}\setminus B_{l_k}$. Consequently $\frac 1{2\pi}\int_{B_{2l_k}}h_i^k e^{u_i^k}>\sigma_i+\epsilon$ for some $\epsilon>0$ independent of $k$, a contradiction to the definition of $\sigma_i$ in (\ref{12sep2e3}).

Let
$$v_i^k(y)=u_i^k(l_ky)+2\log l_k,\quad i\in I. $$
Then clearly we have
\begin{equation}\label{12may4e3}
\left\{\begin{array}{ll}
\Delta v_i^k(y)+\sum_{j=1}^n  a_{ij}H_j^k(y)e^{v_j^k(y)}=0,\quad |y|\le 1/l_k,\quad i\in I\\
\\
\bar v_i^k(1)\to -\infty,
\end{array}
\right.
\end{equation}
where
$$H_i^k(y)=h_i^k(l_ky),\quad i\in I,\quad |y|\le 1/l_k. $$
The Pohozaev identity we use is
\begin{eqnarray}
&& \sum_i \int_{B_{\sqrt{R_k}}}(x\cdot \nabla  H_i^k)e^{v_i^k}+2\sum_i \int_{B_{\sqrt{R_k}}} H_i^k e^{v_i^k} \nonumber \\
&=& \sqrt{R_k}\int_{\partial B_{\sqrt{R_k}}} \sum_i H_i^k e^{v_i^k}+\sqrt{R_k}\int_{\partial B_{\sqrt{R_k}}}\sum_{i,j}\big ( a^{ij}\partial_{\nu}v_i^k\partial_{\nu} v_j^k-\frac 12 a^{ij}\nabla v_i^k \nabla v_j^k\big )
\label{poho1}
\end{eqnarray}
where $R_k\to \infty$ will be chosen later, $( a^{ij})$ is the inverse matrix of $( a_{ij})$. The key point of the following proof is to
choose $R_k$ properly in order to estimate $\nabla v_i^k$ on $\partial B_{\sqrt{R_k}}$. In the estimate of $\partial B_{\sqrt{R_k}}$,
the procedure is to get rid of not important parts and prove that the radial part of $\nabla v_i^k$ is the leading term.
To estimate all the terms of the Pohozaev identity we first write (\ref{poho1}) as
$$\mathcal{L}_1+\mathcal{L}_2=\mathcal{R}_1+\mathcal{R}_2+\mathcal{R}_3 $$
where $\mathcal{L}_1$ stands for ``the first term on the left". Other terms are understood similarly. First we choose $R_k\to \infty$ such that
$R_k^{3/2}=o(l_k^{-1})$, then by using
$l_k\to 0$ to show that $\mathcal{L}_1=o(1)$. To evaluate $\mathcal{L}_2$, we observe that by Lemma \ref{oscillation},
$v_i^k(y)\to -\infty $
over all compact subsets of $\mathbb R^2\setminus B_{1/2}$. Thus we further require $R_k$ to satisfy
\begin{equation}\label{12may24e1}
\int_{B_{R_k}\setminus B_{3/4}} H_i^k e^{v_i^k}=o(1)
\end{equation}
and for $i\in I$, by (\ref{12may4e3}) and Lemma \ref{oscillation}
\begin{equation}\label{12may23e2}
v_i^k(y)+2\log |y|\to -\infty, \mbox{ uniformly in }\quad 1<|y|\le R_k.
\end{equation}
By the choice of $l_k$ we clearly have
$$\frac 1{2\pi}\int_{B_1}H_i^ke^{v_i^k}=\frac 1{2\pi }\int_{B_{l_k}}h_i^k e^{u_i^k}=\sigma_i+o(1),\quad i\in I. $$
By (\ref{12may24e1}) we have
$$\mathcal{L}_2=4\pi \sum_{i=1}^n \sigma_i+o(1). $$
For $\mathcal{R}_1$ we use (\ref{12may23e2}) to have $\mathcal{R}_1=o(1)$.

Therefore we are left with the estimates of $\mathcal{R}_2$ and $\mathcal{R}_3$, for which we shall estimate $\nabla v_i^k$ on $\partial B_{R_k}$. Let
$$G_k(y,\eta)=-\frac 1{2\pi }\log |y-\eta |+\gamma_k(y,\eta) $$
be the Green's function on $B_{l_k^{-1}}$ with respect to Dirichlet boundary condition. Clearly
$$\gamma_k(y,\eta)=\frac 1{2\pi}\log \frac{|y|}{l_k^{-1}}|\frac{l_k^{-2}y}{|y|^2}-\eta | $$
and we have
\begin{equation}\label{12may4e5}
\nabla_y \gamma_k(y,\eta)=O(l_k), \quad y\in \partial B_{\sqrt{R_k}},\quad \eta\in B_{l_k^{-1}}.
\end{equation}

We first estimate $\nabla v_i^k$ on $\partial B_{R_k^{1/2}}$. By Green's representation formula
$$v_i^k(y)=\int_{B_{l_k^{-1}}}G(y,\eta)\sum_{j=1}^n a_{ij}H_i^ke^{v_j^k}d\eta+H_{ik}, $$
where $H_{ik}$ is the harmonic function satisfying $H_{ik}=v_i^k$ on $\partial B_{l_k^{-1}}$. Since $H_{ik}-c_k=O(1)$ for some $c_k$, $|\nabla H_{ik}(y)|=O(l_k),$
\begin{eqnarray}\label{12may4e6}
\nabla v_i^k(y)&=&\int_{B_{l_k^{-1}}}\nabla_y G_k(y,\eta)\sum_{j=1}^n  a_{ij} H_j^k e^{v_j^k}d\eta
+\nabla H_{ik}(y)\\
&=&-\frac 1{2\pi}\int_{B_{l_k^{-1}}}\frac{y-\eta}{|y-\eta |^2}\sum_{j=1}^n  a_{ij}  H_j^k e^{v_j^k}d\eta
+O(l_k).\nonumber
\end{eqnarray}

We estimate the integral in (\ref{12may4e6}) over a few subregions. First the integral over $B_{l_k^{-1}}\setminus B_{R_k^{2/3}}$ is $o(1)R_k^{-\frac 12}$ because over this region $1/|y-\eta|\sim 1/|\eta |\le o(R_k^{-1/2})$.
For the integral over $B_1$, we use
$$\frac{y-\eta}{|y-\eta |^2}=\frac{y}{|y|^2}+O(1/|y|^2)$$
to obtain
$$-\frac 1{2\pi}\int_{B_1}\frac{y-\eta}{|y-\eta |^2}\sum_{j=1}^n  a_{ij} H_j^k e^{v_j^k}
=(-\frac{y}{|y|^2}+O(1/|y|^2))(\sum_{j=1}^n  a_{ij}\sigma_j+o(1)). $$
This is the leading term. For the integral over region $B(0, \sqrt{R_k}/2)\setminus B_1$, we use $1/|y-\eta |\sim 1/|y|$ and (\ref{12may24e1}) to get
$$\int_{B_{R_k^{1/2}/2}\setminus B_1}\frac{y-\eta}{|y-\eta |^2}\sum_{j=1}^n  a_{ij}  H_j^k e^{v_j^k}=o(1)|y|^{-1}. $$
By similar argument we also have
$$\int_{B_{R_k^{2/3}}\setminus (B_{R_k^{1/2}/2}\cup B(y, \frac{|y|}2))}\frac{y-\eta}{|y-\eta |^2}\sum_{j=1}^n  a_{ij}  H_j^k e^{v_j^k}=o(1)|y|^{-1}. $$
Finally over the region $B(y, \frac{|y|}2)$ we use $e^{v_i^k(\eta)}=o(1)|\eta |^{-2}$ to get
$$\int_{B(y, \frac{|y|}2)} \frac{y-\eta }{|y-\eta |^2} \sum_{j=1}^n  a_{ij}  H_j^k e^{v_j^k}=o(1)|y|^{-1}. $$

Combining the estimates on all the subregions mentioned above we have
$$\nabla v_i^k(y)=(-\frac{y}{|y|^2})(\sum_{j=1}^n  a_{ij} \sigma_j+o(1))+o(|y|^{-1}), \quad |y|=R_k^{\frac 12}. $$
Using the above in $\mathcal{R}_2$ and $\mathcal{R}_3$ we have
$$\sum_{i,j=1}^n  a_{ij}\sigma_i\sigma_j=4\sum_{i=1}^n \sigma_i+\circ(1). $$
Proposition \ref{piforU} is established for the non-singular case.

\medskip

{\bf Case two: Singular case: $\exists \gamma_i\neq 0$. }

\medskip

\begin{lem}\label{singpi}
For $\sigma\in (0,1)$, the following Pohozaev identity holds:
\begin{eqnarray*}
&&\sigma\int_{\partial B_{\sigma}}\sum_{i,j\in I}a^{ij}\big (\partial_{\nu}u_i^k\partial_{\nu}u_j^k-\frac{1}2\nabla u_i^k\cdot
\nabla u_j^k\big )+\sum_{i\in I}\sigma\int_{\partial B_{\sigma}}h_i^ke^{u_i^k}\\
&=&2\sum_{i\in I}\int_{B_{\sigma}}h_i^ke^{u_i^k}+\sum_{i\in I}\int_{B_{\sigma}}(x\cdot\nabla h_i^k)e^{u_i^k}+4\pi\sum_{i,j\in I}a^{ij}\gamma_i^k\gamma_j^k.
\end{eqnarray*}
\end{lem}

\noindent{\bf Proof of Lemma \ref{singpi}:}

First, we claim that for each fixed $k$,
\begin{equation}\label{12dec17e1}
\nabla u_i^k(x)=2\gamma_i^kx/|x|^2+O(1)\quad \mbox{ near the origin.}
\end{equation}
 Indeed, recall the equation for the regular part $\tilde u^k_i$ is
$$\Delta \tilde u_i^k(x)+\sum_{j}|x|^{2\gamma_j^k}h_j^k(x)e^{\tilde u_j^k(x)}=0\quad B_1. $$
By the argument of Lemma 4.1 in \cite{linzhang1}, for fixed $k$, $\tilde u_i^k$ is bounded above near $0$, then elliptic estimate leads to (\ref{12dec17e1}).

Let $\Omega=B_{\sigma}\setminus B_{\epsilon}$. Then standard Pohozaev identity on $\Omega$ is
\begin{eqnarray*}
&&\sum_{i\in I}\bigg (\int_{\Omega}(x\cdot \nabla
h_i^k)e^{u_i^k}+2h_i^ke^{u_i^k}\bigg )\\
&=&\int_{\partial \Omega}\bigg (\sum_i(x\cdot
\nu)h_i^ke^{u_i^k}+\sum_{i,j}a^{ij}(\partial_{\nu}u_j^k(x\cdot \nabla
u_i^k)-\frac 12(x\cdot \nu)(\nabla u_i^k\cdot \nabla u_j^k))\bigg
).
\end{eqnarray*}
Let $\epsilon\to 0$, then the integration over $\Omega$ extends to $B_{\sigma}$ by the integrability of $h_i^ke^{u_i^k}$ and (\ref{ah}).
For the terms on the right hand side, clearly $\partial \Omega=\partial B_{\sigma}\cup \partial B_{\epsilon}$. Thanks to (\ref{12dec17e1}), the integral on $\partial B_{\epsilon}$ is $-4\pi \sum_{i,j}a^{ij}\gamma_i^k\gamma_j^k$.
Lemma \ref{singpi} is established. $\Box$

\medskip

Let
$$\sigma_i^k(r)=\frac 1{2\pi }\int_{B_r}h_i^k e^{u_i^k},\quad i\in I, $$
then we have
\begin{lem}\label{spi2}
Let $\epsilon_k\to 0$ such that $\Sigma_k\subset B(0, \epsilon_k/2)$ and
\begin{equation}\label{121217e1}
u_i^k(x)+2\log |x|\to -\infty,\quad |x|=\epsilon_k, \quad i\in I.
\end{equation}
Then we have
\begin{equation}\label{singpi3}
\sum_{i,j\in I}a_{ij}\sigma_i^k(\epsilon_k)\sigma_j^k(\epsilon_k)=4\sum_{i\in I}(1+\gamma_i^k)\sigma_i^k(\epsilon_k)+o(1).
\end{equation}
\end{lem}

\noindent{\bf Proof of Lemma \ref{spi2}:}
First the existence of $\epsilon_k$ that satisfies (\ref{121217e1}) is guaranteed by Lemma \ref{oscillation}. In $B_{\epsilon_k}$, we let
$\tilde u_i^k(x)$ be defined as in (\ref{reg-u}).
Then
$$v_i^k(y)=\tilde u_i^k(\epsilon_ky)+2(1+\gamma_i^k)\log \epsilon_k. $$
Using $v_i^k\to -\infty$ on $\partial B_1$, we obtain, by Green's representation formula and standard estimates,
$$\nabla v_i^k(y)=(\sum_{j\in I}a_{ij}\sigma_j^k(\epsilon_k)+o(1))y, \quad y\in \partial B_1. $$
After translating the above to estimates of $u_i^k$, we have
\begin{equation}\label{121217e2}
\nabla u_i^k(x)=(\sum_{j\in I}a_{ij}\sigma_j^k(\epsilon_k)-2\gamma_j^k)x/|x|^2+o(1)/|x|, \quad |x|=\epsilon_k.
\end{equation}
As we observe the Pohozaev identity in Lemma \ref{singpi} with $\sigma=\epsilon_k$, we see easily that the second term on the LHS and the second term on the RHS are both $o(1)$. The first term on the RHS is clearly $4\pi \sum_i \sigma_i^k(\epsilon_k)$. Therefore we only need to evaluate the first term on the LHS, for which we use (\ref{121217e2}).
Lemma \ref{spi2} is established by similar estimates as in the nonsingular case. $\Box$

\medskip
Proposition \ref{piforU} is established for the singular case as well. $\Box$

\medskip

\begin{rem}\label{rem1}
The proof of Proposition \ref{piforU} clearly indicates the following statements when it is applied to $SU(3)$ Toda system.
Let $B(p_k,l_k)$ be a circle centered at $p_k$ with radius $l_k$. Let $\Sigma_k'$ be a subset of $\Sigma_k$. Suppose
$dist(\Sigma_k',\partial B(p_k,l_k))=o(1)dist(\Sigma_k\setminus \Sigma_k',\partial B(p_k,l_k))$ and we consider the following two situations: If $p_k=0$, we have
$$\tilde \sigma_1^k(l_k)^2-\tilde \sigma_1^k(l_k)\tilde \sigma_2^k(l_k)^2+\tilde \sigma_2^k(l_k)=2\mu_1\tilde \sigma_1^k(l_k)+2\mu_2\tilde \sigma_2^k(l_k)+o(1). $$
If $0\in \Sigma_k\setminus \Sigma_k'$, then
$$\tilde \sigma_1^k(l_k)^2-\tilde \sigma_1^k(l_k)\tilde \sigma_2^k(l_k)+\tilde \sigma_2^k(l_k)^2=2\tilde \sigma_1^k(l_k)+2\tilde \sigma_2^k(l_k)+o(1) $$
where $\tilde \sigma_i^k(l_k)=\frac 1{2\pi}\int_{B(p_k,l_k)}h_i^ke^{u_i^k}$. This fact will be used in the final step of the proof of Theorem \ref{cla-ener}.
\end{rem}

\begin{rem} From the proof of Proposition \ref{piforU} we see that the Pohozaev identity has to be evaluated on \emph{fast decay} components in order to rule out the $\mathcal{R}_1$ term. A component is called \emph{fast decay} if the difference between itself the thresh-hold harmonic function tends to $-\infty$, for example, see (\ref{121217e1}). A component is called a \emph{slow decay} component if it is not a fast decay component. Later in the remaining part of the proof of Theorem \ref{cla-ener} we shall derive Pohozaev identities over different regions and all of them will have to be evaluated on fast decay components.
\end{rem}

\section{Fully bubbling systems}

Next we consider a typical blowup situation for systems: Fully bubbling solutions. First let $\gamma_i^k\equiv 0$ for all $i\in I$. Let
\begin{equation}\label{121223e4}
\lambda^k=\max\{\max_{B_1}u_1^k,...,\max_{B_1}u_n^k\}
\end{equation}
and $x^k\to 0$ be where $\lambda^k$ is attained.  Let
\begin{equation}\label{121223e5}
v_i^k(y)=u_i^k(x_k+e^{-\frac 12\lambda^k}y)-\lambda^k, \quad i\in I, \quad y\in \Omega_k
\end{equation}
where $\Omega_k=\{y;\quad e^{-\frac 12\lambda^k}y+x_k\in B_1 \}$. The sequence is called fully bubbling if, along a subsequence
\begin{equation}\label{121223e1}
\{v_1^k,....,v_n^k\} \mbox{ converge in $C^2_{loc}(\mathbb R^2)$ to $(v_1,...,v_n)$}
\end{equation}
that satisfies
\begin{equation}\label{121223e2}
\Delta v_i+\sum_{j\in I}a_{ij}h_je^{v_j}=0, \quad \mathbb R^2, \quad i\in I
\end{equation}
where $h_i=\lim_{k\to \infty}h_i^k(0)$.
Our next theorem is concerned with the closeness between $u^k=(u_1^k,...,u_n^k)$ and $v=(v_1,..,v_n)$.
\begin{thm}\label{thm3}
Let $A=A_n$, $u^k$ be a sequence of solutions to (\ref{12jun22e8}) with $\gamma_i^k=0, \forall i\in I$. Suppose $u^k$ satisfies
(\ref{12jun27e1}) and (\ref{osc-enr}), $h^k$ satisfies (\ref{ah}), $\lambda^k$, $x^k$, $v^k$ are described by (\ref{121223e4}), (\ref{121223e5}), respectively. Suppose $u^k$ is fully bubbling, then
there exists $C>0$ independent of $k$ such that
\begin{equation}\label{12jun27e2}
|u_i^k(e^{-\frac 12\lambda^k}y+x^k)-\lambda^k-v_i(y)|\le C+o(1)\log (1+|y|), \quad \mbox{ for }x\in \Omega_k,\, i\in I.
\end{equation}
\end{thm}

\begin{rem}
If $A$ is nonnegative, i.e. the system is Liouville system, Theorem \ref{thm3} and Theorem \ref{thm2} below are established in \cite{linzhang1}.
For $A=A_2$, Jost-Lin-Wang \cite{jostlinwang}
proved
$$|u_i^k(e^{-\frac 12\lambda^k}y+x^k)-\lambda^k-v_i(y)|\le C, \quad \mbox{ for }x\in \Omega_k,\, i=1,2. $$
Clearly this estimate is slightly stronger than (\ref{12jun27e2}) for $n=2$. The proof of Jost-Lin-Wang is involved with holonomy theory but the proof of Theorem \ref{thm3} is a simply application of the Pohozaev identity proved in section three.
\end{rem}

\medskip

If $\exists \gamma_i\neq 0$, we let
$$\tilde \lambda^k=\max\{\frac{\max_{B_1} \tilde u_1^k}{(1+\gamma_1^k)},..,\frac{\max_{B_1} \tilde u_n^k}{(1+\gamma_n^k)}\}, $$
and
$$\tilde v_i^k(y)=\tilde u_i^k(e^{-\frac 12\tilde \lambda^k}y)-(1+\gamma_i^k)\tilde \lambda^k $$
for $i\in I$ and $y\in \Omega_k:=\{y;\,\, e^{-\frac 12\tilde \lambda^k}y\in B_1\}. $
We assume
\begin{equation}\label{12jun22e9}
\{\tilde v_1^k,...,\tilde v_n^k\} \mbox{ converge in $C^2_{loc}(\mathbb R^2)$ to $(\tilde v_1,...,\tilde v_n)$}
\end{equation}
that satisfies
\begin{equation}\label{12dec27e1}
\Delta \tilde v_i+\sum_{j=1}^n a_{ij}|x|^{2\gamma_j} h_j e^{\tilde v_j}=0\quad \mathbb R^2, \quad i\in I
\end{equation}
where $h_i=\lim_{k\to \infty}h_i^k(0)$.

\begin{thm}\label{thm2}
Let $A=A_n$, $\tilde u^k$, $\tilde v^k$, $(\tilde v_1,..,\tilde v_n)$, $\tilde\lambda^k$, $\epsilon_k$ and $\Omega_k$ be described as above, $h_i^k$ and $\gamma_i^k$ satisfy (\ref{ah}), then under assumption (\ref{12jun22e9}) there exists $C>0$ independent of $k$ such that
\begin{equation}\label{12jun27e12}
|\tilde u_i^k(e^{-\frac 12\tilde\lambda^k}y)-(1+\gamma_i^k)\tilde\lambda^k-\tilde v_i(y)|\le C+o(1)\log (1+|y|), \quad \mbox{ for }x\in \Omega_k.
\end{equation}
\end{thm}

\noindent{\bf Proof of Theorem \ref{thm3}:}

Recall that $\sigma_i$ is defined in (\ref{12sep2e3}).
By Proposition \ref{piforU} we have
\begin{equation}\label{12may24e4}
\sum_{i,j\in I} a_{ij}\sigma_i\sigma_j=4\sum_{i\in I} \sigma_i.
\end{equation}
On the other hand, let
$$\sigma_{iv}:=\frac 1{2\pi}\int_{\mathbb R^2}h_ie^{v_i}, \quad \mbox{ for }i=1,..,n $$
where $v=(v_1,..,v_n)$ is the limit of the fully bubbling sequence after scaling. Clearly
$\sigma_v=(\sigma_{1v},...,\sigma_{nv})$ also satisfies (\ref{12may24e4}).
We claim that
\begin{equation}\label{12may24e6}
\sigma_i=\sigma_{iv},\quad \mbox{for }i=1,..,n.
\end{equation}

Let $s_i=\sigma_i-\sigma_{vi}$, we obviously have $s_i\ge 0$. The difference between $\sigma$ and $\sigma_v$ on (\ref{12may24e4}) gives
\begin{equation}\label{12aug20e2}
\sum_{i,j\in I}a_{ij}s_is_j+2\sum_{i\in I}(\sum_{j\in I}a_{ij}\sigma_{vj}-2)s_i=0.
\end{equation}
First by Proposition \ref{selectprop} we have $\sum_{j\in I}a_{ij}\sigma_{vj}-2>0$. Next if either $A$ is nonnegative ( $a_{ij}\ge 0$ for all $i,j=1,..,n$) or $A$ is positive definite, we have $\sum_{i,j\in I}a_{ij}s_is_j\ge 0$.
Then (\ref{12aug20e2}) and $s_i\ge 0$ imply (\ref{12may24e6}).

From the convergence from $v_i^k$ to $v_i$ in $C^2_{loc}(\mathbb R^2)$ we can find $R_k\to \infty$ such that
 $$|v_i^k(y)-v_i(y)|=o(1),\quad |y|\le R_k. $$
 For $|y|>R_k$, let
 $$\bar v_i^k(r)=\frac 1{2\pi r}\int_{\partial B_r}v_i^k(y)dS_y. $$
 Then
 $$\frac{d}{dr}\bar v_i^k(r)=\frac{1}{2\pi r}\int_{B_r}\Delta v_i^k=-\frac{1}{2\pi r}\int_{B_r}\sum_{j\in I}a_{ij}h_j^ke^{v_j^k}=-\frac{\sum_ja_{ij}\sigma_j+o(1)}{r}. $$
 Hence
 $$\bar v_i^k(r)=-(\sum_{j\in I}a_{ij}\sigma_j+o(1))\log r+O(1), \quad \mbox{for all } r>2. $$
 Since
 $v_i^k(y)=\bar v_i^k(|y|)+O(1)$ and
 $$v_i(y)=-(\sum_j a_{ij}\sigma_j)\log |y|+O(1)\quad \mbox{ for }\quad |y|>1, $$
  we see that (\ref{12jun27e2}) holds.  Theorem \ref{thm3}
 is established. $\Box$

\medskip

\noindent{\bf Proof of Theorem \ref{thm2}:}  By (\ref{singpi3}) we have
\begin{equation}\label{12dec27e2}
\sum_{i,j\in I}a_{ij}\sigma_i\sigma_j=4\sum_{i\in I}(1+\gamma_i)\sigma_i.
\end{equation}
Recall that $v=(v_1,...,v_n)$ satisfies (\ref{12dec27e1}). Let
$$\sigma_{iv}=\frac 1{2\pi}\int_{\mathbb R^2}h_i|x|^{2\gamma_i}e^{v_i}.$$
On one hand, $(\sigma_{1v},...,\sigma_{iv})$ also satisfies (\ref{12dec27e2}), on the other hand, the classification theorem of Lin-Wei-Ye \cite{lin-wei-ye} gives
\begin{equation}\label{12dec27e3}
\sum_{j\in I}a_{ij}\sigma_{jv}>2+2\gamma_i,\quad i\in I.
\end{equation}
Let $s_i=\sigma_i-\sigma_{iv}$ ($i\in I$), then (\ref{12dec27e2}), which is satisfied by both $(\sigma_1,..,\sigma_n)$ and $(\sigma_{1v},...,\sigma_{nv})$, gives
$$\sum_{i,j\in I}a_{ij}s_is_j+2\sum_{i\in I}(\sum_{j\in J}a_{ij}\sigma_{jv}-2-2\gamma_i)s_i=0. $$
By (\ref{12dec27e3}) and the assumption on $A$, we have $s_i=0$ for all $i\in I$. The remaining part of the proof is exactly like the
last part of the proof of Theorem \ref{thm3}. Theorem \ref{thm2} is established. $\Box$

\medskip

\section{Asymptotic behavior of solutions in each simple blowup area}
In this section we derive some results on the energy classification around each blowup point. First we let $A=A_n$ (the Cartan Matrix) and consider

\begin{center}
{\bf The neighborhood around $0$.}
\end{center}

 Since $0$ is postulated to belong to $\Sigma_k$ first, it means there may not be a bubbling picture in a neighborhood of $0$.

Let $\tau_k=\frac 12 dist(0, \Sigma_k\setminus \{0\})$ we consider the energy limits of $h_i^ke^{u_i^k}$ in $B_{\tau_k}$. By the selection process and Lemma \ref{oscillation},
\begin{equation}\label{si1}
u_i^k(x)+2\log |x|\le C, \quad u_i^k(x)=\bar u_i^k(|x|)+O(1) \quad i\in I, \quad |x|\le \tau_k
\end{equation}
where $\bar u_i^k(|x|)$ is the average of $u_i^k$ on $\partial B_{|x|}$.
Let $\tilde u_i^k$ be defined by (\ref{reg-u}).
Then we have
$$\Delta \tilde u_i^k(x)+\sum_{j\in I}a_{ij}|x|^{2\gamma_j}h_j^k(x)e^{\tilde u_j^k(x)}=0, \quad |x|\le \tau_k. $$

Let
$$-2\log \delta_k=\max_{i\in I}\max_{x\in B(0,\tau_k)}\frac{\tilde u_i^k}{1+\gamma_i^k}$$
and
\begin{equation}\label{13jan16e1}
v_i^k(y)=\tilde u_i^k(\delta_ky)+2(1+\gamma_i^k)\log \delta_k, \quad |y|\le \tau_k/\delta_k.
\end{equation}

East to see the equation for $v_i^k$ is
$$\Delta v_i^k(y)+\sum_{j\in I}a_{ij}|y|^{2\gamma_j^k}h_j^k(\delta_ky)e^{v_j^k(y)}=0, \quad |y|\le \tau_k/\delta_k. $$

Then we consider two trivial cases. First: $\tau_k/\delta_k\le C$. This is a case that there is no entire bubble after scaling.

Let $f_i^k$ solve
$$\left\{\begin{array}{ll}
\Delta f_i^k+\sum_{j\in I}a_{ij}|y|^{2\gamma_j^k}h_j^k(\delta_ky)e^{v_j^k}=0, \quad |y|\le \tau_k/\delta_k, \\
f_i^k=0, \quad \mbox{ on }\quad |y|=\tau_k/\delta_k.
\end{array}
\right.
$$
Using $v_i\le 0$ we have $|f_i^k|\le C$ on $B(0, \tau_k/\delta_k)$. Since $v_i^k-f_i^k$ is harmonic and $v_i^k$ has bounded oscillation on
$\partial B(0, \tau_k/\delta_k)$, we have
\begin{equation}\label{si2}
v_i^k(x)=\bar v_i^k(\partial B(0, \tau_k/\delta_k))+O(1), \quad \forall x\in B(0, \tau_k/\delta_k)
\end{equation}
where $\bar v_i^k(\partial B(0, \tau_k/\delta_k))$ stands for the average of $v_i^k$ on $\partial B(0, \tau_k/\delta_k)$.
Direct computation shows that
$$\int_{B(0,\tau_k)}e^{u_i^k(x)}dx=\int_{B(0,\tau_k/\delta_k)}e^{v_i^k(y)}|y|^{2\gamma_i^k}dy. $$
Therefore
\begin{equation}\label{13jan2e1}
\int_{B_{\tau_k}}h_i^ke^{u_i^k}dx=O(1)e^{\bar v_i^k(\partial B(0, \tau_k/\delta_k))}.
\end{equation}
So if $\bar v_i^k(\partial B(0, \tau_k/\delta_k))\to -\infty$, $\int_{B_{\tau_k}}h_i^ke^{u_i^k}dx=o(1)$.

\medskip

The second trivial case is when the blowup sequence is fully bubbling.  Clearly we now have
\begin{equation}\label{13jan2e2}
\tau_k/\delta_k\to \infty
\end{equation}
and we assume that
$(v_1^k,...,v_n^k)\to (v_1,...,v_n)$ in $C^2_{loc}(\mathbb R^2)$. Clearly
$$\Delta v_i+\sum_{j=1}^n a_{ij}|x|^{2\gamma_j}h_j e^{v_j}=0\quad \mathbb R^2, \quad i\in I $$
where $h_i=\lim_{k\to \infty}h_i^k(0)$.
By the classification theorem of Lin-Wei-Ye \cite{lin-wei-ye}, we have
$$\frac 1{2\pi }\sum_{j\in I}a_{ij}\int_{\mathbb R^2}|y|^{2\gamma_j}e^{v_j}h_jdy=2(2+\gamma_i+\gamma_{n+1-i}) $$
and
$$v_i(y)=-(4+2\gamma_{n+1-i})\log |y| +O(1), \quad |y|>1,\quad i\in I. $$
By the proof of Theorem \ref{thm2} that there is only one bubble.

\medskip

The final case we consider is a partially blown-up picture. Note that (\ref{13jan2e2}) is assumed. For the following two propositions we assume $n=2$. i.e. we consider $SU(3)$ Toda systems.

\begin{prop}\label{singat0} Suppose (\ref{12jun22e8}), (\ref{12jun27e1}), (\ref{ah}) and (\ref{osc-enr}) hold for $u^k$, $h_i^k$ and $\gamma_i$ etc. The matrix $A=A_2$. (\ref{13jan2e2}) also holds.
Suppose $s_k\in (0,\tau_k)$ satisfies
$$u_i^k(x)\le -2\log |x|-N_k, \quad i=1,2$$
for all $|x|=s_k$ and some $N_k\to \infty$. Then $(\sigma_1^k(s_k),\sigma_2^k(s_k))$ is an $o(1)$ perturbation of one of the following five types:
\begin{eqnarray*}
&&(2\mu_1,0),\quad (0,2\mu_2),\quad (2(\mu_1+\mu_2),2\mu_2),\\
&&(2\mu_1,2(\mu_1+\mu_2)),\quad
(2\mu_1+2\mu_2,2\mu_1+2\mu_2).
\end{eqnarray*}

On $\partial B(0,\tau_k)$, for each $i$ either
$$u_i^k(x)+2\log |x|\ge -C,\quad |x|=\tau_k$$
for some $C>0$ or
\begin{equation}\label{121223e8}
u_i^k(x)+2\log |x|<-(2+\delta)\log |x|+\delta \log \delta_k,\quad |x|=\tau_k
\end{equation}
for some $\delta>0$.
If (\ref{121223e8}) holds for some $i$, then
$$
\sigma_i^k(\tau_k)=o(1), 2\mu_i+o(1),\mbox{ or } 2\mu_1+2\mu_2+o(1).
$$
Moreover, there exists at least one $i_0$ such that (\ref{121223e8}) holds for $i_0$.
\end{prop}

Similarly for bubbles away from the origin we have

\begin{prop}\label{singaw} Suppose (\ref{12jun22e8}), (\ref{12jun27e1}), (\ref{ah}) and (\ref{osc-enr}) hold for $u^k$, $h_i^k$ and $\gamma_i$ etc. The matrix $A=A_2$.
Let $x_k\in \Sigma_k\setminus \{0\}$,$\bar \tau_k=\frac 12dist(x_k,\Sigma_k\setminus \{0, x_k\})$ and
$$\bar \delta_k=exp\bigg (-\frac 12\max_{i=1,2;x\in B(x_k,\bar \tau_k)}u_i^k(x)\bigg ). $$
Then for all $s_k\in (0, \bar \tau_k)$, if
$$u_i^k(x)+2\log |x-x_k|\le -N_k, \quad \forall i\in I,\quad |x-x_k|=s_k, \quad i=1,2$$
for some $N_k\to \infty$, then $(\frac 1{2\pi}\int_{B(x_k,s_k)}h_1^ke^{u_1^k},\frac 1{2\pi}\int_{B(x_k,s_k)}h_2^ke^{u_2^k})$ is an $o(1)$ perturbation of one of the following five types:
$$(2,0),(0,2),(2,4),(4,2),(4,4). $$
On $\partial B(x_k,\bar \tau_k)$, for each $i$ either
$$u_i^k(x)+2\log \bar \tau_k\ge -C, \quad \forall x\in \partial B(x_k,\bar \tau_k)$$
or
\begin{equation}\label{121223e9}
u_i^k(x)\le -(2+\delta)\log \bar \tau_k+\delta \log \bar \delta_k,\quad \forall x\in \partial B(x_k,\bar \tau_k).
\end{equation}
If (\ref{121223e9}) holds for some $i$, then $\frac{1}{2\pi}\int_{B(x_k,\bar \tau_k)}h_i^ke^{u_i^k}$ is $o(1)$,$2+o(1)$ or $4+o(1)$.
Moreover, there exists at least one $i_0$ such that (\ref{121223e9}) holds for $i_0$.
\end{prop}

We shall only prove Proposition \ref{singat0} as the proof for Proposition \ref{singaw} is similar.

\noindent{\bf Proof of Proposition \ref{singat0}:}

Let $v_i^k$ be defined by (\ref{13jan16e1}). Since we only need to consider a partially blown-up situation, without loss of generality we assume $v_1^k$ converges to $v_1$ in $C^2_{loc}(\mathbb R^2)$ and $v_2^k$ tends to $-\infty$ over any compact subset of $\mathbb R^2$. The equation for $v_1$ is
$$\Delta v_1+2h_1|y|^{2\gamma_1}e^{v_1}=0,\quad \mathbb R^2, \quad \int_{\mathbb R^2}h_1|y|^{2\gamma_1}e^{v_1}<\infty. $$
where $h_1=\lim_{k\to \infty} h_1^k(0)$.
By the classification result of Prajapat-Tarantello \cite{prajapat} we have
$$2\int_{\mathbb R^2} h_1|y|^{2\gamma_1}e^{v_1}=8\pi \mu_1$$
and
$$v_1(y)=-4\mu_1\log |y|+O(1),\quad |y|>1. $$
Thus we can find $R_k\to \infty$ (without loss of generality $R_k=o(1)\tau_k/\delta_k$) such that
$$\frac 1{2\pi}\int_{B_{R_k}}h_1^k(\delta_ky)|y|^{2\gamma_1^k}e^{v_1^k}=2\mu_1+o(1),\quad  (i.e. \quad \sigma_1^k(\delta_k R_k)=2\mu_1+o(1)) $$
and
$$\int_{B_{R_k}}h_2^k(\delta_ky)|y|^{2\gamma_2^k}e^{v_2^k}=o(1). $$
For $r\ge R_k$, recall that,
$$\sigma_i^k(\delta_kr)=\frac 1{2\pi}\int_{B_r}h_i^k(\delta_ky)|y|^{2\gamma_i^k}e^{v_i^k}dy $$
then we have
\begin{eqnarray*}
\frac{d}{dr}\bar v_1^k(r)&=&\frac{-2\sigma_1^k(\delta_k r)+\sigma_2^k(\delta_k r)}r, \\
\frac{d}{dr}\bar v_2^k(r)&=&\frac{\sigma_1^k(\delta_k r)-2\sigma_2^k(\delta_k r)}r \quad R_k\le r\le \tau_k/\delta_k.
\end{eqnarray*}
Clearly we have
\begin{equation}\label{13jan26e1}
R_k\frac{d}{dr}\bar v_1^k(R_k)=-4\mu_1+o(1), \quad R_k\frac{d}{dr}\bar v_2^k(R_k)=2\mu_1+o(1).
\end{equation}

The following lemma says that as long as both components stay well below the harmonic function $-2\log |y|$ (i.e. both of them are fast decay components), there is no essential change on the energy for either component:
\begin{lem}\label{ener-cn}
Suppose $L_k\in (R_k,\tau_k/ \delta_k)$ satisfies
\begin{equation}\label{13jan24e1}
v_i^k(y)+2\gamma_i^k\log |y|\le -2\log |y|-N_k, \quad R_k\le |y|\le L_k, \quad i=1,2
\end{equation}
for some $N_k\to \infty$, then
$$\sigma_{i}^k(\delta_kR_k)=\sigma_{i}^k(\delta_kL_k)+o(1), \quad i=1,2. $$
\end{lem}

\noindent{\bf Proof of Lemma \ref{ener-cn}:} We aim to prove that $\sigma_{i}^k$ does not change much from $\delta_kR_k$ to $\delta_kL_k$. Suppose this is not the case, then there exists $i$ such that $\sigma_{i}^k(\delta_kL_k)>\sigma_{i}^k(\delta_kR_k)+\delta$ for some $\delta>0$. Let $\tilde L_k\in (R_k, L_k)$ such that
\begin{equation}\label{13mar12e8}
\max_{i=1,2}(\sigma_{i}^k(\delta_k\tilde L_k)-\sigma_{i}^k(\delta_kR_k))=\epsilon, \quad \forall i=1,2
\end{equation}
where $\epsilon>0$ is sufficiently small. Then for $v_1^k$,
\begin{equation}\label{13mar10e2}
\frac{d}{dr}\bar v_1^k(r)\le \frac{-4(1+\gamma_1)+\epsilon}r\le -\frac{2(1+\gamma_1)+\epsilon}r.
\end{equation}
Then it is easy to see from Lemma \ref{oscillation} that
$$\int_{B_{\tilde L_k}\setminus B_{R_k}}|y|^{2\gamma_1^k}e^{v_1^k}=o(1), $$
which is $\sigma_1^k(\delta_k \tilde L_k)=\sigma_1^k(\delta_k R_k)+o(1)$. Indeed, by Lemma \ref{oscillation}
$$\int_{B_{L_k}\setminus B_{R_k}}|y|^{2\gamma_1^k}e^{v_1^k}=O(1)\int_{B_{L_k}\setminus B_{R_k}}|y|^{2\gamma_1^k}e^{\bar v_1^k}=o(1). $$
The second equality above is because by (\ref{13mar10e2})
$$\bar v_1^k(r)+2\gamma_1^k\log r\le -N_k-2\log R_k+(-2-\epsilon/2)\log r,\quad R_k\le r\le L_k. $$

 Thus $\sigma_2^k(\delta_k \tilde L_k)=\sigma_2^k(\delta_kR_k)+\epsilon$.
However, since (\ref{13jan24e1}) holds, by Remark \ref{rem1} we have
$$\lim_{k\to \infty}(\sigma_1^k(\delta_k\tilde L_k),\sigma_2^k(\delta_k\tilde L_k))\in \Gamma. $$
The two points on $\Gamma$ that have the first component equal to $2\mu_1$ are $(2\mu_1,0)$ and $(2\mu_1,2(\mu_1+\mu_2))$. Thus (\ref{13mar12e8})
is impossible.  Lemma \ref{ener-cn} is established. $\Box$

\medskip

From Lemma \ref{ener-cn} and (\ref{13jan26e1}) we see that for $r\ge R_k$, either
\begin{equation}\label{13jan24e3}
v_i^k(y)+2\gamma_i^k\log |y|\le -2\log |y|-N_k, \quad R_k\le |y|\le \tau_k/\delta_k, \quad i=1,2
\end{equation}
or there exists $L_k\in (R_k,\tau_k/\delta_k)$ such that
\begin{equation}\label{13jan15e1}
v_{2}^k(y)+2\gamma_2^k\log L_k\ge -2\log L_k-C \quad |y|=L_k
\end{equation}
for some $C>0$, while for $R_k\le |y|\le L_k$,
\begin{equation}\label{13jan24e2}
v_1^k(y)+2\gamma_1^k\log |y|\le -(2+\delta)\log |y|, \quad R_k\le |y|\le L_k
\end{equation}
for some $\delta>0$.
Indeed, from (\ref{13jan26e1}) we see that if the energy has to change, $\sigma_2^k$ has to change first. $L_k$ can be chosen so that
$\sigma_2^k(\delta_kL_k)-\sigma_2^k(\delta_kR_k)=\epsilon$ for some $\epsilon>0$ small.

\begin{lem}\label{13lem1} Suppose
there exist $L_k\ge R_k$  such that (\ref{13jan15e1}) and (\ref{13jan24e2}) hold. For $L_k$ we assume $L_k=o(1)\tau_k/\delta_k$.
Then there exists $\tilde L_k$ such that $\tilde L_k/L_k\to \infty$ and $\tilde L_k=o(1)\tau_k/\delta_k$ still holds. For $|y|=\tilde L_k$, we have
\begin{equation}\label{12dec31e1}
v_i^k(y)+2(1+\gamma_i^k)\log |y|\le -N_k,\quad |y|=\tilde L_k,\quad i=1,2
\end{equation}
for some $N_k\to \infty$.
In particular
\begin{equation}\label{12dec31e2}
v_1^k(y)+2(1+\gamma_1^k+\frac{\delta}4)\log |y|\le 0,\quad |y|=\tilde L_k.
\end{equation}
\begin{equation}\label{13jan1e2}
\sigma_1^k(\delta_k\tilde L_k)=2\mu_1+o(1),\quad \sigma_2^k(\delta_k\tilde L_k)=2\mu_1+2\mu_2+o(1).
\end{equation}
\end{lem}

\begin{rem} The statement of Lemma \ref{13lem1} can be understood as follows: Suppose starting from $\partial B_{L_k}$, $\sigma_2^k$ starts to change because (\ref{13jan15e1}) holds. Then from $L_k$ to $\tilde L_k$, $\sigma_1^k$ does not change much and $v_1^k$ is still way below $-2(1+\gamma_1^k)\log |y|$ but $v_2^k$ has changed from decaying slowly (which is (\ref{13jan15e1})) to a fast decay ( the $i=2$ part of (\ref{12dec31e2})). In other words, as $\sigma_2^k$ changes from $L_k$ to $\tilde L_k$, $v_2^k$ changes from slow decay to fast decay but $v_1^k$ still has fast decay in the meanwhile. The change of $\sigma_2^k$ has influenced the derivative of $\bar v_1^k$ but has not made $\sigma_1^k$ change much because $\sigma_2^k$ changes too fast from $L_k$ to $\tilde L_k$.
\end{rem}

\noindent{\bf Proof of Lemma \ref{13lem1}:} First we observe that by Lemma \ref{ener-cn} the energy does not change if both components satisfy
(\ref{13jan24e3}). Thus we can assume that $\sigma_2^k(\delta_kL_k)\le \epsilon$ for some $\epsilon>0$ small. Consequently
$$\frac{d}{dr}\bar v_1^k(r)\le \frac{-4(1+\gamma_1)+2\epsilon}r,\quad r\ge R_k. $$
Now we claim that there exists $N>1$ such that
\begin{equation}
\label{sigma2k}
\sigma_2^k(\delta_k(L_k N))\ge 2+\gamma_1+\gamma_2+o(1).
\end{equation}
If this is not true, we would have $\epsilon_0>0$ and $\tilde R_k\to \infty$ such that
\begin{equation}\label{13jan24e5}
\sigma_2^k(\delta_k \tilde R_k L_k)\le 2+\gamma_1+\gamma_2-\epsilon_0.
\end{equation}
On the other hand $\tilde R_k$ can be chosen to tend to infinity slowly so that, by Lemma \ref{oscillation} and (\ref{13jan24e2})
\begin{equation}\label{13jan26e2}
v_1^k(y)+2(1+\gamma_1^k)\log |y|\le -\frac{\delta}2 \log |y|, \quad L_k\le |y|\le \tilde R_k L_k.
\end{equation}
Clearly (\ref{13jan26e2}) implies $\sigma_1^k(\delta_kL_k)=\sigma_1^k(\delta_k\tilde R_kL_k)+o(1)$. Thus by (\ref{13jan24e5})
\begin{equation}\label{13jan24e4}
\frac{d}{dr}\bar v_2^k(r)\ge \frac{-2-2\gamma_2+\epsilon_0/2}r.
\end{equation}
Using (\ref{13jan24e4}) and
$$v_2^k(y)=(-2-2\gamma_2^k)\log |y|+O(1), \quad |y|=L_k $$
we see easily that
$$\int_{B(0, \tilde R_kL_k)\setminus B(0, L_k)}|y|^{2\gamma_2^k}e^{v_2^k}\to \infty, $$
a contradiction to (\ref{osc-enr}).Therefore (\ref{sigma2k}) holds.

By Lemma \ref{oscillation}
$$v_i^k(y)+2\log (NL_k)=\bar v_i^k(NL_k)+2\log (NL_k)+O(1), \quad i=1,2,\quad |y|=NL_k. $$
Thus we have
\begin{eqnarray*}
&& \bar v_1^k(NL_k)\le (-2-2\gamma_1^k-\delta/2)\log (NL_k), \\
&& \bar v_2^k(NL_k)\ge (-2-2\gamma_2^k)\log (NL_k)-C.
\end{eqnarray*}
Consequently
$$\bar v_2^k((N+1)L_k)\ge (-2-2\gamma_2^k)\log L_k-C, $$
leads to
$$\frac 1{2\pi}\int_{B(0, (N+1)L_k)}h_2^k(\delta_ky)|y|^{2\gamma_2^k}e^{v_2^k(y)}dy\ge 2+\gamma_1+\gamma_2+\epsilon_0 $$
for some $\epsilon_0>0$. Going back to the equation for $\bar v_2^k$ we have
$$\frac{d}{dr}\bar v_2^k(r)\le -\frac{2+2\gamma_2+\epsilon_0}r, \quad r= (N+1)L_k. $$
Therefore we can find $\tilde R_k\to \infty$ such that $\tilde R_kL_k=o(1)\tau_k/\delta_k$ and
\begin{eqnarray*}
v_2^k(y)\le (-2-2\gamma_2^k-\epsilon_0)\log |y|-N_k,\quad |y|=\tilde R_k L_k, \\
v_1^k(y)\le (-2-2\gamma_1^k-\delta/4)\log |y|,\quad L_k\le |y|\le \tilde R_kL_k.
\end{eqnarray*}
Obviously
$$\sigma_1^k(\delta_k\tilde R_kL_k)=\sigma_1^k(\delta_k L_k)+o(1)=\sigma_1^k(\delta_k R_k)+o(1)=2(1+\gamma_1)+o(1). $$
By computing the Pohozaev identity on $\tilde R_k L_k$ we have
$$\sigma_2^k(\delta_k \tilde R_k L_k)=2\mu_1+2\mu_2+o(1). $$
Letting $\tilde L_k=\tilde R_kL_k$ we have proved Lemma \ref{13lem1}. $\Box$

\medskip

To finish the proof of Proposition \ref{singat0} we need to consider the region $\tilde L_k\le |y|\le \tau_k/\delta_k$ if $L_k=o(1)\tau_k/\delta_k$ (in which case $\tilde L_k$ can be made as $o(1)\tau_k/\delta_k$), or $L_k=O(1)\tau_k/\delta_k$.  First we consider the region
$\tilde L_k\le |y|\le \tau_k/\delta_k$ when $\tilde L_k=o(1)\tau_k/\delta_k$.
It is easy to verify that
\begin{eqnarray*}
\frac{d}{dr}\bar v_1^k(r)=-\frac{2\gamma_1-2\gamma_2}r+o(1)/r,\quad r=\tilde L_k, \\
\frac{d}{dr}\bar v_2^k(r)=-\frac{6+2\gamma_1+4\gamma_2+o(1)}r,\quad r=\tilde L_k.
\end{eqnarray*}
The second equation above implies
$$\frac{d}{dr}\bar v_2^k(r)\le -\frac{2\mu_2+\delta}r, \quad r=\tilde L_k $$
for some $\delta>0$. So $\sigma_2^k(r)$ does not change for $r\ge \tilde L_k$ unless $\sigma_1^k$ changes. By the same argument as before, either $v_1^k$ rises to $-2\log |y|+O(1)$ on $|y|=\tau_k/\delta_k$ or there is $\hat L_k=o(1)\tau_k/\delta_k$ such that
$$\sigma_i^k(\delta_k\hat L_k)=2\mu_1+2\mu_2+o(1), \quad i=1,2. $$
Since this is the energy of a fully blowup system, we have in this case both
$$v_i^k(y)\le -(2\mu_i+\delta)\log |y|, \quad |y|=\tau_k\delta_k ,\quad i=1,2$$
for some $\delta>0$.

If $L_k=O(1)\tau_k/\delta_k$. In this case it is easy to use Lemma \ref{oscillation} to see that one component is $-2(1+\gamma_i^k)\log |y|+O(1)$ and the other component has the fast decay. Proposition \ref{singat0} is established. $\Box$

\section{Combination of bubbling areas}

The following definition plays an important role:

\begin{Def} \label{group1} Let $Q_k=\{p_1^k,..,p_q^k\}$ be a subset of $\Sigma_k$ such that $Q_k$ has more than one point in it and $\Sigma_k\setminus Q_k
=\not \emptyset$. $Q_k$ is called a group if and
\begin{enumerate}
\item
$$dist(p_i^k,p_j^k)\sim dist(p_s^k,p_t^k), $$
where $p_i^k,p_j^k,p_s^k,p_t^k$ are any points in $Q_k$ such that $p_i^k\neq p_j^k$ and $p_t^k\neq p_s^k$.
\item For any $p_k\in \Sigma_k\setminus Q_k$,
$\frac{dist(p_i^k,p_j^k)}{dist(p_i^k,p_k)}\to 0$ for all $p_i^k,p_j^k\in Q_k$ with $p_i^k\neq p_j^k$.
\end{enumerate}
\end{Def}

{\bf Proof of Theorem \ref{cla-ener}:} Let $2\tau_k$ be the distance between $0$ and $\Sigma_k\setminus \{0\}$. For each $z_k\in \Sigma_k\cap
\partial B(0,2\tau_k)$, if $dist(z_k,\Sigma_k\setminus \{z_k\})\sim \tau_k$, let $G_0$ be the group that contains the origin. On the other hand, if there exists $z_k'\in \partial B(0,2\tau_k)$ such that $\tau_k/dist(z_k',\Sigma_k\setminus z_k')\to \infty$ we let $G_0$ be $0$ itself. By the definition of group, all members of $G_0$ are in $B(0,N\tau_k)$ for some $N$ independent of $k$. Let
$$ v_i^k(y)=u_i^k(\tau_ky)+2\log \tau_k, \quad |y|\le \tau_k^{-1}. $$
Then we have
\begin{equation}\label{13jan17e3}
\Delta  v_i^k(y)+\sum_{j=1}^2a_{ij}h_j^k(\tau_ky)e^{v_j^k(y)}=4\pi \gamma_i^k\delta_0,\quad |y|\le \tau_k^{-1}.
\end{equation}

Let $0$, $Q_1$,...,$Q_m$ be the images of members of $G_0$ after the scaling from $y$ to $\tau_ky$. Then all $Q_i\in B_N$. By Proposition \ref{singat0} and Proposition \ref{singaw} at least one component decays fast on $\partial B_1$. Without loss of generality we assume
$$v_1^k\le -N_k\quad \mbox{ on }\quad \partial B_1 $$
for some $N_k\to \infty$ and
$$\sigma_1^k(\tau_k)=o(1),2\mu_1+o(1)\mbox{ or } 2\mu_1+2\mu_2+o(1). $$
Specifically, if $\tau_k/\delta_k\le C$, $\sigma_1^k(\tau_k)=o(1)$. Otherwise, $\sigma_1^k(\tau_k)$ is equal to the two other cases mentioned above.
By Lemma \ref{oscillation} $v_1^k\le -N_k+C$ on all $\partial B(Q_t,1)$ ($t=1,...,m$), therefore by Proposition \ref{singaw},
$$\frac 1{2\pi}\int_{B(Q_t,1)}h_1^k(\tau_k\cdot)e^{v_1^k}=2m_t+o(1), \quad t=1,...,m $$
where for each $t$, $m_t=0,1$ or $2$.  Let $2\tau_kL_k$ be the distance from $0$ to the nearest group other than $G_0$. Then $L_k\to \infty$. By Lemma \ref{oscillation} and the proof of Lemma \ref{12aug8lem1} we can find $\tilde L_k\le L_k$, $\tilde L_k\to \infty$ such that most of the energy of $v_1^k$ in $B(0,\tilde L_k)$ is contributed by bubbles and $v_2^k$ decays faster than $-2\log \tilde L_k$ on $\partial B(0,\tilde L_k)$:
\begin{eqnarray}\label{13mar15e1}
&&\frac 1{2\pi}\int_{B(0,\bar L_k)}h_1^k(0)e^{v_1^k}\\
&=&2m+o(1), \quad 2\mu_1+2m+o(1)\quad \mbox{ or }\quad  2(\mu_1+\mu_2)+2m+o(1)\nonumber
\end{eqnarray}
for some nonnegative integer $m$, and
\begin{equation}\label{13mar11e1}
v_2^k(y)+2\log \tilde L_k\to -\infty \quad |y|=\tilde L_k.
\end{equation}

Then we evaluate the Pohozaev identity on $B(0, \tilde L_k)$. Since (\ref{13mar11e1}) holds, by Remark \ref{rem1} we have
$$\lim_{k\to \infty} (\sigma_1^k(\tau_k\tilde L_k),\sigma_2^k(\tau_k\tilde L_k))\in \Gamma. $$
Moreover, by (\ref{13mar15e1})
we see that $\lim_{k\to \infty} (\sigma_1^k(\tau_k\tilde L_k),\sigma_2^k(\tau_k\tilde L_k))\in \Sigma$ because the limit point is the intersection
between the line $\sigma_1=\lim_{k\to \infty}\sigma_1^k(\tau_k\tilde L_k)$ with $\Gamma$.

The Pohozaev identity for
$(\sigma_1^k(\tau_k\tilde L_k),\sigma_2^k(\tau_k\tilde L_k))$ can be written as
\begin{eqnarray*}
&&\sigma_1^k(\tau_k\tilde L_k)(2\sigma_1^k(\tau_k\tilde L_k)-\sigma_2^k(\tau_k\tilde L_k)-4\mu_1)\\
&&+\sigma_2^k(\tau_k\tilde L_k)(2\sigma_2^k(\tau_k\tilde L_k)-\sigma_1^k(\tau_k\tilde L_k)-4\mu_2)=o(1).
\end{eqnarray*}
Thus either
\begin{equation}\label{13jan27e1}
2\sigma_1^k(\tau_k\tilde L_k)-\sigma_2^k(\tau_k\tilde L_k)\ge 4\mu_1+o(1)
\end{equation}
or
$$2\sigma_2^k(\tau_k\tilde L_k)-\sigma_1^k(\tau_k\tilde L_k)\ge 4\mu_2+o(1). $$
Moreover,
if
$$2\sigma_1^k(\tau_k\tilde L_k)-\sigma_2^k(\tau_k\tilde L_k)\ge 2\mu_1+o(1) \mbox{ and }
2\sigma_2^k(\tau_k\tilde L_k)-\sigma_1^k(\tau_k\tilde L_k)\ge 2\mu_2+o(1), $$
by the proof of Theorem \ref{thm2}
$$\int_{B_{l_k}\setminus \tau_k\tilde l_k}h_i^ke^{u_i^k}=o(1),\quad i=1,2$$
for any $l_k\to 0$. In this case we have
$$\sigma_i=\lim_{k\to \infty}\sigma_i^k(\tau_k\tilde L_k),\quad i=1,2$$
and Theorem \ref{cla-ener} is proved in this case.

Thus without loss of generality we assume that (\ref{13jan27e1}) holds. From the equation for $u_1^k$, this means for some $\delta>0$
\begin{equation}\label{13jan27e2}
\bar u_1^k(\tau_k\tilde L_k)\le -2\log (\tau_k\tilde L_k)-N_k,\quad
\frac{d}{dr}\bar u_1^k(r)<(-2-\delta)/r,\quad r=\tau_k\tilde L_k.
\end{equation}
The property above implies, by the proof of Proposition \ref{singat0}, that as $r$ grows from $\tau_k\tilde L_k$ to $\tau_kL_k$,
the following three situations may occur:

{\bf Case one: } Both $u_i^k$ satisfy, for some $N_k\to \infty$, that
$$u_i^k(x)+2\log |x|\le -N_k, \quad \tau_k\tilde L_k\le |x|\le \tau_kL_k, \quad i=1,2. $$
In this case
$$\sigma_i^k(\tau_k\tilde L_k)=\sigma_i^k(\tau_k L_k)+o(1), \quad i=1.2 $$
So on $\partial B(0,\tau_kL_k)$, $u_1^k$ is still a fast decaying component.

{\bf Case two:} There exist $L_{1,k}, L_{2,k}\in (\tilde L_k, L_k)$ such that
$$u_2^k(x)\ge -2\log (\tau_kL_{1,k})-C \quad |x|=\tau_kL_{1,k},$$
\begin{equation}\label{13mar12e1}
u_i^k(x)\le -2\log (\tau_kL_{2,k})-N_k\quad |x|=\tau_K L_{2,K}, \quad i=1,2
\end{equation}
and
\begin{equation}\label{13mar12e2}
\sigma_1^k(\tau_k \tilde L_{k})=\sigma_1^k(\tau_k L_{2,k})+o(1).
\end{equation}
Since (\ref{13mar12e1}) holds, by Remark \ref{rem1},
$(\lim_{k\to \infty}\sigma_1^k(\tau_kL_{2,k}),\lim_{k\to \infty}\sigma_2^k(\tau_k L_{2,k}))\in \Gamma$. Then we further observe that since (\ref{13mar12e2}) holds, $\lim_{k\to \infty}(\sigma_1^k(\tau_kL_{2,k}),\sigma_2^k(\tau_k L_{2,k}))\in \Sigma$ because this point is obtained by intersecting $\Gamma$ with $\sigma_1=\lim_{k\to \infty} \sigma_1^k(\tau_k\tilde L_k)$. In other words, the new point
$\lim_{k\to \infty}(\sigma_1^k(\tau_kL_{2,k}),\sigma_2^k(\tau_k L_{2,k}))$ is on the upper right part of the old point
$\lim_{k\to \infty}(\sigma_1^k(\tau_k\tilde L_k),\sigma_2^k(\tau_k \tilde L_k))$.

{\bf Case three:}
$$u_2^k(x)\ge -2\log \tau_k L_k-C, \quad |x|=\tau_k L_k $$
for some $C>0$ and $\sigma_1^k(\tau_k \tilde L_k)=\sigma_1^k(\tau_k L_k)+o(1)$. This means at $\partial B(0,\tau_k L_k)$, $u_1^k$ is still the fast decaying component.

\medskip

If the second case above happens, the discussion of the relationship between $\sigma_1^k$ and $\sigma_2^k$ on $B(0,\tau_k L_k)\setminus
B(0,\tau_k L_{2,k})$ is the same as before.
In any case on $\partial B(0,\tau_kL_k)$ at least one of the two components has fast decay and has its energy equal to a corresponding component of a point in $\Sigma$. For any group not equal to $G_0$, it is easy to see that the fast decay component has its energy equal to $0$, $2$ or $4$. The combination of bubbles for groups is very similar to the combination of bubbling disks as we have done before. For example, let $G_0,G_1,...,G_t$ be
groups in $B(0,\epsilon_k)$ for some $\epsilon_k\to 0$. Suppose the distance between any two of $G_0,..,G_t$ are comparable and
$$dist(G_i,G_j)=o(1)\epsilon_k, \quad \forall i,j=0,.., t, \quad i\neq j. $$
Also we require $(\Sigma_k\setminus (\cup_{i=0}^tG_i)) \cap B(0,2\epsilon_k)=\emptyset$. Let $\epsilon_{1,k}=dist(G_0,G_1)$, then all $G_0,...,G_t$
are in $B(0,N\epsilon_{1,k})$ for some $N>0$. Without loss of generality let $u_1^k$ be a fast decaying component on $\partial B(0,N\epsilon_{1,k})$. Then we have
$$\sigma_1^k(N\epsilon_{1,k})=\sigma_1^k(\tau_k L_k)+2m+o(1) $$
where $m$ is a nonnegative integer because by Lemma \ref{oscillation}, $u_1^k$ is also a fast decaying component for $G_1,...,G_t$. Moreover, by Proposition \ref{singaw}, the energy of $u_1^k$ in $G_s$ ($s=1,...,t$) is $o(1),2+o(1)$ or $4+o(1)$. If $u_2^k$ also has fast decay on
$\partial B(0, N\epsilon_{1,k})$, then $\lim_{k\to \infty}(\sigma_1^k(N\epsilon_{1,k}), \sigma_1^k(N\epsilon_{1,k}))\in \Sigma$ because this is a point of intersection between $\Gamma$ and $\sigma_1=\lim_{k\to \infty}\sigma_1^k(\tau_k L_k)+2m$.
If
$$u_2^k(x)\ge -2\log N \epsilon_{1,k}-C,\quad |x|=N\epsilon_{1,k}, $$
then as before we can find $\epsilon_{3,k}$ in $(N\epsilon_{1,k},\epsilon_k)$ such that, for some $N_k\to \infty$,
$$u_i^k(x)+2\log \epsilon_{3,k}\le -N_k,\quad i=1,2, \quad |x|=\epsilon_{3,k} $$
and
$$\sigma_1^k(N\epsilon_{1,k})=\sigma_1^k(\epsilon_{3,k}). $$
Thus we have
$$\lim_{k\to \infty}(\sigma_1^k(\epsilon_{3,k}),\sigma_2^k(\epsilon_{3,k}))\in \Sigma. $$
because this point is the intersection
between $\Gamma$ and $\sigma_1=\lim_{k\to \infty}\sigma_1^k(N\epsilon_{1,k})$.

The last possibility on $B(0,\epsilon_k)\setminus B(0,\epsilon_{1,k})$ is
$$\sigma_1^k(\epsilon_{k})=\sigma_1^k(N\epsilon_{1,k})+o(1)$$ and
$$u_2^k(x)+2\log \epsilon_{k}\ge -C, \quad |x|=\epsilon_{k}. $$
In this case $u_1^k$ is the fast decaying component on $\partial B(0, \epsilon_k)$.

Such a procedure can be applied to include groups further away from $G_0$. Since we have only finite blowup disks this procedure only needs to be applied finite times. Finally let $s_k\to 0$ such that
$$\sigma_i=\lim_{k\to \infty}\lim_{s_k\to 0}\sigma_i^k(s_k),\quad i=1,2 $$
and, for some $N_k\to \infty$,
$$u_i^k(x)+2\log s_k\le -N_k,\quad i=1,2,\quad |x|=s_k. $$
Then we see that $(\sigma_1,\sigma_2)\in \Sigma$.
Theorem \ref{cla-ener} is established. $\Box$


\begin{thebibliography}{99}

\bibitem{bclt} Bartolucci, D.; Chen, Chiun-Chuan; Lin, Chang-Shou; Tarantello, Gabriella Profile of blow-up solutions to mean field equations with singular data.  {\em Comm. Partial Differential Equations}  29  (2004),  no. 7-8, 1241-1265.
\bibitem{bart} Bartolucci, D.; Lin, Chang-shou;
Sharp existence results for mean field equations with singular data.
{\em J. Differential Equations} 252 (2012), no. 7, 4115–4137.
\bibitem{bart2} Bartolucci, D.; Tarantello, G.
The Liouville equation with singular data: a concentration-compactness principle via a local representation formula.
{\em J. Differential Equations} 185 (2002), no. 1, 161–180.
\bibitem{bart3} Bartolucci, D.; Tarantello, G.
Liouville type equations with singular data and their applications to periodic multivortices for the electroweak theory.
{\em Comm. Math. Phys. } 229 (2002), no. 1, 3–47.
\bibitem{bart4}  Bartolucci, D.; Malchiodi, A. 
 An improved geometric inequality via vanishing moments, with applications to singular Liouville equations. {\em Comm. Math. Phys.}  322  (2013),  no. 2, 415–452. 
\bibitem{malchiodi-b} Battaglia, L; Malchiodi, A; A Moser-Trudinger Inequality for the singular Toda system, {\em preprint}.
\bibitem{bennet} W. H. Bennet, Magnetically self-focusing streams, {\em Phys. Rev.} 45 (1934), 890-897.
\bibitem{bolton1} Bolton, J., Woodward, L.M.: Some geometrical aspects of the 2-dimensional Toda equations.
In: Geometry, Topology and Physics, Campinas, 1996, pp. 69–81. de Gruyter, Berlin
(1997).
\bibitem{bolton2} Bolton, J., Jensen, G.R., Rigoli, M.,Woodward, L.M.: On conformal minimal immersions
of S2 into CPn. {\em Math. Ann. } 279(4), 599–-620 (1988).
\bibitem{calabi} Calabi, E.: Isometric imbedding of complex manifolds. {\em Ann. Math.} 58(2), 1–-23 (1953).
\bibitem{chanillo} S. Chanillo, M. K-H Kiessling, Conformally invariant systems of nonlinear PDE of
Liouville type. {\em Geom. Funct. Anal. } 5 (1995), no. 6, 924-947.
\bibitem{chen-lin-2} C. C. Chen, C. S. Lin, Estimate of the conformal scalar curvature equation via the method of moving planes. II.
 {\em J. Differential Geom. } 49 (1998), no. 1, 115–178.
\bibitem{chenlin1} C. C. Chen, C. S. Lin, Sharp estimates for solutions of multi-bubbles in compact Riemann surfaces.  {\em Comm. Pure Appl. Math. }  55  (2002),  no. 6, 728-771.
\bibitem{chenlin2} C. C. Chen, C. S. Lin, Topological degree for a mean field equation on Riemann surfaces.
{\em Comm. Pure Appl. Math. } 56 (2003), no. 12, 1667–1727.
\bibitem{lin} C. C. Chen, C. S. Lin, Estimate of the conformal scalar curvature equation via the method of moving planes. II. {\em  J. Differential Geom.} 49 (1998), no. 1, 115–178.
\bibitem{chenliduke} W. X. Chen, C. M. Li, Classification of solutions of some nonlinear elliptic equations. {\em Duke Math. J. } 63 (1991), no. 3, 615–-622.
\bibitem{chen-li-duke2} W. X. Chen, C. M. Li, Qualitative properties of solutions to some nonlinear elliptic
equations in $R\sp 2$. {\em Duke Math. J. } 71 (1993), no. 2, 427--439.
\bibitem{chern} Chern, S.S., Wolfson, J.G.: Harmonic maps of the two-sphere into a complex Grassmann
manifold. II. {\em Ann. Math.} 125(2), 301–-335 (1987).
\bibitem{child} S. Childress and J. K. Percus, Nonlinear aspects of Chemotaxis, {\em Math. Biosci. }    56
(1981), 217--237.
\bibitem{chipot} M. Chipot, I. Shafrir, G. Wolansky, On the solutions of Liouville systems. {\em J. Differential
Equations} 140 (1997), no. 1, 59--105.
\bibitem{debye} P. Debye and E. Huckel, Zur Theorie der Electrolyte, {\em Phys. Zft} 24 (1923), 305--325.
\bibitem{doliwa} Doliwa, A.: Holomorphic curves and Toda systems. {\em Lett. Math. Phys.} 39(1), 21–-32 (1997).
\bibitem{dunne1} G. Dunne, R. Jackie, S.Y.Pi, C. Trugenberger, Self-dual Chern-Simons solitons and two dimensional nonlinear equations, {\em Phys. Rev.
D} 43 (1991), 1332--1345.
\bibitem{dunne2} G. Dunne, Self-dual Chern-Simons theories. Lecture Notes in Physics, Springer, Berline, 1995.
\bibitem{ganoulis} N. Ganoulis, P. Goddard, D.  Olive, : Self-dual monopoles and Toda molecules. {\em   Nucl.
Phys. B } 205, 601–-636 (1982)
\bibitem{guest} Guest, M.A.: Harmonic Maps, Loop Groups, and Integrable Systems. London Mathematical
Society Student Texts, vol. 38. Cambridge University Press, Cambridge (1997).
\bibitem{jostlinwang} J. Jost, C. S. Lin and G. F. Wang, Analytic aspects of the Toda system II: bubbling behavior and existence of solutions,
{\em Comm. Pure Appl. Math.} 59 (2006), no. 4, 526--558.
\bibitem{jostwang} J. Jost, G. F. Wang, Classification of solutions of a Toda system in R2. {\em Int. Math. Res. Not.} 2002, no. 6, 277–-290.
\bibitem{keller} E. F. Keller and L. A. Segel, Traveling bands of Chemotactic Bacteria: A theoretical
analysis, {\em J. Theor. Biol.} 30 (1971), 235--248.
\bibitem{kiessling} M. K.-H. Kiessling and J. L. Lebowitz, Dissipative stationary Plasmas: Kinetic Modeling Bennet Pinch, and generalizations, {\em Phys. Plasmas} 1 (1994), 1841--1849.
\bibitem{leznov1} Leznov, A.N.: On the complete integrability of a nonlinear system of partial differential
equations in two-dimensional space. {\em Theor. Math. Phys.} 42, 225–-229 (1980).
\bibitem{leznov2} Leznov, A.N., Saveliev, M.V.: Group-Theoretical Methods for Integration of Nonlinear
Dynamical Systems. Progress in Physics, vol. 15. Birkh\"auser, Basel (1992)
\bibitem{licmp} Y. Y. Li, Harnack type inequality: the method of moving planes,
{\em Comm. Math. Phys.} 200 (1999), no. 2, 421--444.
\bibitem{li} Y. Y. Li, Prescribing scalar curvature on Sn and related problems. I.{\em J. Differential Equations } 120 (1995), no. 2, 319–410.
\bibitem{lin-wei-ye} C.S. Lin, J. C. Wei, D. Ye, Classifcation and nondegeneracy of $SU(n+1)$ Toda system, {\em Invent. Math. } 190(2012), no.1, 169-207.
\bibitem{lin-wei-zhao} C. S. Lin, J. C. Wei, C. Zhao, Sharp estimates for fully bubbling solutions of a $SU(3)$ Toda system, {\em Geom. Funct. Anal. } 22 (2012), no. 6, 1591–1635.
\bibitem{linzhang1} C. S. Lin, L. Zhang, Profile of bubbling solutions to a Liouville system , {\em Annales de l'Institut Henri Poincare / Analyse non lineaire} Volume 27, Issue 1, January-February 2010, Pages 117--143,
\bibitem{linzhang2} C. S. Lin, L. Zhang, A topological degree counting for some Liouville systems of mean field type.
{\em Comm. Pure Appl. Math. } 64 (2011), no. 4, 556–-590.
\bibitem{linzhang3} C. S. Lin, L. Zhang, On Liouville systems at critical parameters, Part 1: One bubble. {\em J. Funct. Anal.} 264 (2013), no. 11, 2584–2636.
\bibitem{linzhang4} C. S. Lin, L. Zhang, On Liouville systems at critical parameters, Part 2: Multi bubbles, in preparation.
\bibitem{MN} A. Malchiodi, C.B. Ndiaye, Some existence results for the
Toda system on closed surfaces, {\em Att. Accad. Naz. Lincei
Cl. Sci. Fis. Mat. Natur. Rend. Lincei (9) Mat. Appl.} 18(2007), no.4,
391-412.
\bibitem{MR} A. Malchiodi and D. Ruiz, a variational analysis of the Toda
system on compact surfaces, {\em Comm. Pure Appl. Math.}  66 (2013), no. 3, 332–371.
\bibitem{mansfield} Mansfield, P.: Solutions of Toda systems. {\em Nucl. Phys. B} 208, 277–-300 (1982).
\bibitem{mock} M. S. Mock, Asymptotic behavior of solutions of transport equations for semiconductor devices, {\em J. Math. Anal. Appl.} 49 (1975), 215--225.
\bibitem{nolasco1} Nolasco, M., Tarantello, G.: Double vortex condensates in the Chern-Simons theory. {\em Calc.
Var. Partial Differ. Equ.} 9, 31–-94 (1999).
\bibitem{nolasco2} Nolasco, M., Tarantello, G.: Vortex condensates for the SU(3) Chern-Simons theory. {\em Commun.
Math. Phys. } 213(3), 599–-639 (2000).
\bibitem{pmw}A. Pistoia, M. Musso and J. Wei, New concentration phenomena for $SU(3)$ Toda system, {\em preprint}.

\bibitem{prajapat} J. Prajapat, G. Tarantello, On a class of elliptic problems in $\mathbb R^2$: symmetry and uniqueness results. {\em Proc. Roy. Soc. Edinburgh Sect. A} 131 (2001),no. 4, 967--985.
\bibitem{schoen}. R. Schoen, Stanford lecture notes.

\bibitem{wzz} J. C. Wei, C. Y. Zhao, F. Zhou, On nondegeneracy of solutions of $SU(3)$ Toda system {\em  CRAS } 349(2011), no.3-4, 185--190.
\bibitem{yang1} Yang, Y.: The relativistic non-abelian Chern-Simons equation. {\em Commun. Math. Phys.} 186(1),
199–-218 (1999)
\bibitem{yang2} Yang, Y.: Solitons in Field Theory and Nonlinear Analysis. SpringerMonographs in Mathematics.
Springer, New York (2001)
\bibitem{zhangcmp} L. Zhang, Blowup solutions of some nonlinear elliptic equations
involving exponential nonlinearities.
{\em Comm. Math. Phys.} 268 (2006), no. 1, 105--133.
\bibitem{zhangccm} L. Zhang, Asymptotic behavior of blowup solutions for elliptic equations with exponential nonlinearity and singular data. {\em Commun. Contemp. Math. } 11  (2009),  no. 3, 395--411.

\end{thebibliography}
\end{document}